\documentclass{article}

\usepackage{amsmath,amsfonts,amssymb,url,fullpage,theorem,mathtools}

\usepackage[utf8]{inputenc}

{\theorembodyfont{\slshape}
\newtheorem{proposition}{Proposition}[section]
\newtheorem{lemma}[proposition]{Lemma}
\newtheorem{corollary}[proposition]{Corollary}
\newtheorem{theorem}[proposition]{Theorem}}

{\theorembodyfont{\upshape}
\newtheorem{remark}[proposition]{Remark}

\newtheorem{example}[proposition]{Example}}

\numberwithin{equation}{section}

\newcommand{\eps}{\varepsilon}
\newcommand{\ph}{\varphi}

\newcommand\C{{\mathbb C}}
\newcommand\poincare{{\mathbb H}}
\renewcommand\P{{\mathbb P}}
\newcommand\Q{{\mathbb Q}}
\newcommand\R{{\mathbb R}}
\newcommand\Z{{\mathbb Z}}

\newcommand{\DD}{{\mathcal D}}

\newcommand\Gal{{\mathrm{Gal}}}
\newcommand\GL{{\mathrm{GL}}}
\newcommand\nlc{{\mathrm{nlc}}}
\newcommand\ord{{\mathrm{ord}}}
\renewcommand\Im{{\mathrm{Im}}}

\newcommand\SL{{\mathrm{SL}}}

\newcommand\tilf{{\widetilde{f}}}
\newcommand\tilg{{\widetilde{g}}}

\newcommand\qed{\hfill$\square$}

1

\makeatletter

\renewcommand*\l@section[2]{%
  \ifnum \c@tocdepth >\z@
    \addpenalty\@secpenalty
    \addvspace{0.2em \@plus\p@}%
    \setlength\@tempdima{1.5em}%
    \begingroup
      \parindent \z@ \rightskip \@pnumwidth
      \parfillskip -\@pnumwidth
      \leavevmode \bfseries
      \advance\leftskip\@tempdima
      \hskip -\leftskip
      #1\nobreak\hfil \nobreak\hb@xt@\@pnumwidth{\hss #2}\par
    \endgroup
  \fi}
  
 \makeatother

\title{Collinear CM-points}

\author{Yuri Bilu, Florian Luca and David Masser}

\date\today

\allowdisplaybreaks[4]

\begin{document}

\hfuzz 3pt

%\vfuzz 100pt

%\hbadness 9999

\maketitle

\begin{abstract}

Andr\'e's celebrated Theorem of 1998 implies that each complex straight line  ${Ax+By+C=0}$ (apart from obvious exceptions) contains at most finitely many points $(j(\tau),j(\tau'))$,  where ${\tau,\tau'\in\poincare}$ are algebraic of degree~$2$. We show that there are only a finite number of such lines which contain more than two such points. As there is a line through any two complex points, this is best possible.

\end{abstract}

{\footnotesize \tableofcontents}

\section{Introduction}
\label{sintr}

In 1998 André~\cite{An98} proved that a non-special irreducible plane curve in $\C^2$ may have at most finitely many CM-points. Here a \textsl{plane curve} is a curve defined by an irreducible  equation ${F(x,y)=0}$, where~$F$ is a polynomial with complex coefficients, and \textsl{CM-point} (called also \textsl{special point}) in $\C^2$ is a point whose coordinates are both singular moduli. Recall that a \textsl{singular modulus} is the invariant of an elliptic curve with complex multiplication; in other words, it is an algebraic number of the form $j(\tau)$, where~$j$ denotes   the standard $j$-function on the upper half-plane~$\poincare$ and ${\tau\in \poincare}$ is an algebraic number of degree~$2$. 
Thus, a $CM$-point is a point   of the form $(j(\tau),j(\tau'))$ with   ${\tau,\tau'\in \poincare}$ algebraic of degree~$2$.

\textsl{Special curves} are those of the following types:

\begin{itemize}
\item
``vertical lines'' ${x=j(\tau)}$ and ``horizontal lines'' ${y=j(\tau)}$, where $j(\tau)$ is a singular modulus;

\item
\textsl{modular curves} $Y_0(N)$, realized as the plane curves ${\Phi_N(x,y)=0}$, where $\Phi_N$ is the modular polynomial of level~$N$. 
\end{itemize}

Recall that the polynomial  ${\Phi_N(X,Y)\in \C[X,Y]}$ is the $X$-monic $\C$-irreducible polynomial satisfying ${\Phi_N(j(z),j(Nz))=0}$. It is known that actually ${\Phi_N(X,Y)\in \Z[X,Y]}$; this and other properties of $\Phi_N$ can be found, for instance, in \cite[Theorem~11.18]{Co89}.

Clearly, each special curve contains infinitely many CM-points, and André proved that special curves are characterized by this property.

André's result was the first non-trivial contribution to the celebrated André-Oort conjecture on the special subvarieties of Shimura varieties; see~\cite{Pi11} and the references therein. 

Several other proofs (some conditional on GRH) of André's theorem were suggested, see \cite{BMZ13,Br01,Ed98,Ku12,Ku13,Pi09}. We mention specially the  argument of Pila~\cite{Pi09},  based on an idea of Pila and Zannier~\cite{PZ08}. In~\cite{Pi11} Pila extended it to   higher dimensions, proving the André-Oort conjecture for subvarieties of $\C^n$. To state this result, one needs to introduce the notion of ``special variety''; then Pila's theorem asserts that an algebraic subvariety of $\C^n$ has at most finitely many maximal special subvarieties. See Section~\ref{spila} and  Theorem~\ref{thpila} for the details. 

Besides general results, some particular curves were considered. For instance, Kühne \cite[Theorem~5]{Ku13} proved that the straight line ${x+y=1}$ has no CM-points\footnote{The same result was  independently obtained in an earlier version of \cite{BMZ13}, but did not appear in the final version.}, and a similar result for the hyperbola ${xy=1}$ was obtained in \cite{BMZ13}. Similarly in~\cite{JHM16} for the quartic curve
{%\small
\begin{gather*}
x^3y-2x^2y^2+xy^3-1728x^3+1216x^2y+1216xy^2-1728y^3+3538944x^2-2752512xy+3538944y^2\\
-2415919104x-2415919104y+549755813888=0;
\end{gather*}
}%
this is equivalent to the fact that there are no complex $t \neq 0,1,-1$ for which the two elliptic curves ${Y^2=X(X-1)(X-t)}$ and ${Y^2=X(X-1)(X+t)}$ both have complex multiplication.

%See also~\cite{BLP14} for a generalization to the infinite family   of curves ${xy=A}$ with ${A\in \Q}$. 

One can ask about CM-points on general straight lines %defined over~$\Q$ (that is, defined by an equation 
${Ax+By+C=0}$. 
%where ${A_1,A_2,B\in \Q}$). 
One has to exclude from consideration the \textsl{special straight lines}: ${x=j(\tau)}$, ${y=j(\tau)}$ (where $j(\tau)$ is a singular modulus) and ${x=y}$, the latter being nothing else than the modular curve $Y_0(1)$ (the modular polynomial $\Phi_1$ is ${X-Y}$).  According to the theorem of André, these are the only straight lines containing infinitely many CM-points.

In~\cite{ABP15} all CM-points lying on non-special straight lines defined over~$\Q$ are listed. 
More generally, Kühne remarks on page~5 of his article~\cite{Ku13} that, given a positive integer~$\nu$, at most finitely many CM-points belong to the union of all non-special straight lines defined over a number field of degree~$\nu$; moreover, for a fixed~$\nu$ all these points can,  in principle, be listed explicitly, though the implied calculation does not seem to be feasible. 

Here we take a different point of view: instead of restricting the degree of field of definition, we study the (non-special) straight lines passing through at least~$3$ CM-points.

Such lines do exist \cite[Remark~5.3]{ABP15}: since 
{%\small
$$
\det\begin{bmatrix*}[l]
1728& -884736000\\
287496& -147197952000 
\end{bmatrix*}=0,
$$}%
the three points $(0,0)$, $(1728,287496)$ and $(-884736000,-147197952000)$ belong to the same straight line ${1331x=8y}$, and just as well for the points $(0,0)$, $(1728,-884736000)$ and $(287496,-147197952000)$ on $512000x=-y$. 
Here
{%\small
\begin{align*}
&j\left(\frac{-1+\sqrt{-3}}{2}\right)=0, \quad j(\sqrt{-1})=1728, \quad j(2\sqrt{-1})=287496,\\
&j\left(\frac{-1+\sqrt{-43}}{2}\right)=-884736000, 
\quad j\left(\frac{-1+\sqrt{-67}}{2}\right)=-147197952000. 
\end{align*}}%

%In the present note we prove that there are only finitely many non-special straight lines passing through~$3$ or more CM-points. We express it in the following equivalent form.

Call an (unordered) triple ${\{P_1,P_2,P_3\}}$ of  CM-points  \textsl{collinear} if  $P_1,P_2,P_3$ are pairwise distinct and  belong to a non-special straight line. %In particular, the points~\eqref{epofl} form a collinear triple; we call the triples~\eqref{epofl} \textsl{exceptional}. 

In this paper we prove the following. 

\begin{theorem}
\label{thcoli}
There exist at most finitely many  collinear triples of CM-points. 
\end{theorem}

In particular, there exist at most finitely many non-special straight lines passing through three or more CM-points. This latter consequence looks formally weaker than Theorem~\ref{thcoli}, but in fact it is equivalent to it, due to the theorem of André.

\begin{remark}
The referee drew our attention to the phenomenon of \textsl{automatic uniformity}, discovered by Scanlon~\cite{Sc04}. Combining Theorem~4.2 from~\cite{Sc04} with Pila's Theorem~\ref{thpila} stated in the next section, one obtains the following ``uniform'' version of the Theorem of André: there is a (non-effective) uniform upper bound~$c_d$ on the number of CM-points
in an arbitrary non-special curve of  geometric degree~$d$ (with an arbitrary
field of definition). For every~$d$, it is a widely open question what the optimal $c_d$
actually is; moreover, even obtaining an effective upper bound for $c_d$ seems to be quite difficult. 
It might be an easier question to ask 
for an optimal bound~$c_d^\ast$ 
such that
\textsl{all but finitely many} non-special curves of degree~$d$ contain at most $c_d^\ast$ 
special points. In this language our Theorem~\ref{thcoli} simply asserts that ${c_1^\ast=2}$. 

\end{remark}

The idea of the proof of Theorem~\ref{thcoli} is simple. Three points $(x_i,y_i)$ lie on a line if and only if 
{%\small
\begin{equation}
\label{edetvar}
\begin{vmatrix}
1&1&1\\
x_1&x_2&x_3\\
y_1&y_2&y_3
\end{vmatrix}=0.
\end{equation}}% 
This defines a variety in $\C^6$  to which we can apply Pila's André-Oort result. This guarantees finiteness outside the special subvarieties of positive dimension. One easily detects ``obvious'' positive-dimensional special subvarieties:  they correspond to the line being special in two dimensions, or the three points not being distinct.  The main difficulty is showing that there are no other positive-dimensional special subvarieties: this is the content of the ``Main Lemma'', whose proof occupies the overwhelming part of the article.  Along the way we have to solve some auxiliary problems not only of André-Oort type but also of ``mixed type'' involving roots of unity.

It could be mentioned that, while the Main Lemma is completely effective, Theorem~\ref{thcoli} is not, because its deduction from the Main Lemma relies on  Pila's Theorem~\ref{thpila}, which is non-effective.  

For analogous Diophantine assertions about lines proved also using ``determinant varieties'', the reader can consult the articles of Evertse, Gy\H{o}ry, Stewart and Tijdeman~\cite{EGST88} about $S$-units, or  of Schlickewei and Wirsing~\cite{SW97} about heights.
In these papers one is actually in the multiplicative group $\mathbb{G}_m^2$ and the appropriate special varieties are much easier to describe.

\paragraph{Plan of the article}
In Section~\ref{spila} we recall the general notion of special variety and state the already mentioned Theorem of Pila, proving the André-Oort conjecture for subvarieties of $\C^n$.

In Section~\ref{sml} we %remark that three points $(x_i,y_i)$ (where ${i=1,2,3}$) belong to a straight line 
state the Main Lemma, which lists all maximal positive-dimensional special subvarieties of the ``determinant variety'' defined by~\eqref{edetvar}, 
and we deduce Theorem~\ref{thcoli} from the Theorem of Pila and the Main Lemma.

In Sections~\ref{sroots},~\ref{smodi},~\ref{srmat} and~\ref{sjmaps} we obtain various auxiliary results used in the sequel. The proof of the Main Lemma occupies Sections~\ref{sini} to~\ref{smggngl}. In Section~\ref{sini} we collect some preliminary material and show how the proof of the Main Lemma splits into four cases. These cases are treated in  Sections~\ref{smg=} to~\ref{smggngl}.

\paragraph{Acknowledgments} Yuri Bilu was supported by  the \textsl{Agence National de la Recherche} project ``Hamot'' (ANR 2010 BLAN-0115-01). We thank Bill Allombert, Qing Liu, Pierre Parent, Jonathan Pila and Thomas Scanlon for useful discussions. 
We also thank the referee, who did the hard job of verifying the proof, detected a number of inaccuracies and made many helpful suggestions.

\section{Special Varieties and the Theorem of Pila}
\label{spila}

%Unless the contrary is stated explicitly,  in this section~$\tau$ with or without indices  denotes an algebraic number of degree~$2$ belonging to~$\poincare$. %, and a capital letter~$N$ with or without indices denotes a positive integer. 

We recall the definition of special varieties from~\cite{Pi11}. The referee pointed out that this is not the definition used in the
standard formulation of the André-Oort conjecture, and some work is required to show that the two are equivalent. However, this presents no issue for our purposes, since the main result that we need, Pila's Theorem~\ref{thpila}, proved in~\cite{Pi11}, is stated therein in terms of this definition. 

To begin with, we define sets~$M$ in $\C^m$ (where ${m \geq 1}$) as follows. If ${m=1}$ then ${M=\C}$, while if ${m \geq 2}$ then~$M$ is given by modular equations
\begin{equation}\label{eqm}\Phi_{N(i)}(x_1,x_i)=0\qquad(i=2,\ldots,m).
\end{equation}

More generally for $\C^n$ (where ${n \geq 1}$) one takes a partition ${n=l_0+m_1+\cdots+m_d}$  (where ${d \geq 0}$) with ${l_0 \geq 0}$ and with  ${m_1 \geq 1,\ldots,m_d \geq 1}$ (when ${d \geq 1}$), and defines sets~$K$ in ${\C^n=\C^{l_0} \times \C^{m_1} \times \cdots \times \C^{m_d}}$ as ${L_0 \times M_1 \times \cdots \times M_d}$, where~$L_0$ (if ${l_0 \geq 1}$) is a single point whose coordinates are singular moduli and ${M_1,\ldots,M_d}$ (if ${d \geq 1}$) are as~$M$ above. Then any irreducible component~$\widetilde K$ of~$K$, which necessarily has the form
\begin{equation}
\label{eqcomp}
\widetilde K=L_0 \times \widetilde M_1 \times \cdots \times \widetilde M_d
\end{equation}
with irreducible components $\widetilde M_1, \ldots, \widetilde M_d$ of $M_1,\ldots,M_d$, is an example of a special variety in the sense of Pila; and one gets all examples by permuting the coordinates. The dimension is $d$.

When ${n=2}$ and ${d=1}$ this agrees with the notion of special curve introduced in Section~\ref{sintr}, because the polynomials $\Phi_N$ are irreducible.

The following property of special varieties  is certainly known, but we could not  find a suitable reference. 

\begin{proposition}
\label{pspec1}
Let $0 \leq e \leq d \leq n$. Then every special variety of dimension $d$ contains a Zariski dense union of special varieties of dimension $e$.
\end{proposition} 

\paragraph{Proof}

If ${d=0}$ there is nothing to prove. Otherwise by induction it suffices to treat the case ${e=d-1}$, with the special variety~\eqref{eqcomp}.

If ${m_1=1}$ then ${\widetilde M_1=\C}$ and for each singular modulus $\xi$ the variety $L_0 \times \{\xi\} \times \widetilde M_2 \times \cdots \times \widetilde M_d$ is special of dimension $d-1$. As there are infinitely many singular moduli, the union is Zariski dense in $\widetilde K$.

If ${m_1 \geq 2}$ (call it $m$) we note from~\eqref{eqm} that~$x_1$ is non-constant on $\widetilde M_1$. Thus the corresponding projection of $\widetilde M_1$ to $\C$ is dominant. We can therefore find infinitely many singular moduli $\xi_1$ for which some $(\xi_1,\xi_2,\ldots,\xi_m)$ lies in $\widetilde M_1$. As ${\Phi_{N(i)}(\xi_1,\xi_i)=0}$ for ${i=2,\ldots,m}$,  it is clear that $\xi_2,\ldots,\xi_m$ are also singular moduli, and now the corresponding
$$L_0 \times \{(\xi_1,\xi_2,\ldots,\xi_m)\} \times \widetilde M_2 \times \cdots \times \widetilde M_d$$
do the trick. \qed

\bigskip

Special points are exactly those of the form ${(\xi_1,\ldots, \xi_n)}$, where each~$\xi_i$ is a singular modulus. To characterize the special curves in a similar way, it will be convenient to use the language of ``$j$-maps''. A map ${f:\poincare\to \C}$ will be called a \textsl{$j$-map} if either ${f(z)=j(\gamma z)}$  
for some ${\gamma \in \GL_2^+(\Q)}$ (a \textsl{non-constant $j$-map}), or ${f(z)=j(\tau)}$ with  ${\tau\in \poincare}$ algebraic of degree~$2$  (a \textsl{constant $j$-map}). Here $\GL_2^+(\Q)$ is the subgroup of $\GL_2(\Q)$ consisting of matrices with positive determinants. We define a \textsl{$j$-set} to be of the form ${\{(f_1(z), \ldots, f_n(z)): z\in \poincare\}}$, where each $f_k$ is a $j$-map and at least one of them is non-constant.%rational $2\times2$ matrices with positive determinant %\footnote{In Section~\ref{sprlm} it will be convenient to attribute level~$0$ to the constant $j$-maps.}).

\begin{remark}
\label{rautom}
It is worth noting that every $j$-map is $\Gamma(N)$-automorphic\footnote{Recall that 
$\Gamma(N)$ is the kernel of the $\bmod N$ reduction map 
${\SL_2(\Z)\to\SL_2(\Z/N\Z)}$, and ``the function~$f$ is $\Gamma(N)$-automorphic'' means ${f\circ \eta=f}$ for any ${\eta \in \Gamma(N)}$.} for some positive integer~$N$. This is trivially true for constant $j$-maps, and a non-constant $j$-map ${f=j\circ\gamma}$ is ${\gamma^{-1}\Gamma(1)\gamma}$-automorphic. So it remains to note that ${\gamma^{-1}\Gamma(1)\gamma}$ contains $\Gamma(N)$ for a suitable~$N$. Indeed, write ${A\in \Gamma(N)}$ as ${I+ NB}$, where~$I$ is the identity matrix and~$B$ is a matrix with entries in~$\Z$. Then the matrix ${\gamma A\gamma^{-1}=I+N\gamma B\gamma^{-1}}$ has entries in~$\Z$ if~$N$ is divisible by the product of the  denominators of the entries of~$\gamma$ and~$\gamma^{-1}$. 
\end{remark}

It seems to be known (and even  used in several places) that every special
curve is a $j$-set and that the converse is also true. As we could not find
a convincing reference, we indicate here an argument. We thank the referee for many explanations on this topic.

\begin{proposition}
\label{pspec2}
\begin{enumerate}
\item
\label{izarclo}
Any $j$-set is a Zariski-closed irreducible algebraic subset of $\C^n$.

\item
\label{ispec2}
A subset of $\C^n$ is a $j$-set if and only if it is a special
curve.
\end{enumerate}

\end{proposition} 

\paragraph{Proof}
In the proof of Part~\ref{izarclo} we may restrict to the case when all $f_1,\ldots, f_n$ are non-constant $j$-maps. Denote by ${Z\subset \C^n}$ the $j$-set defined by these maps.  According to Remark~\ref{rautom}, the maps $f_1,\ldots, f_n$ are $\Gamma(N)$-automorphic for some positive integer~$N$. Hence each~$f_i$ induces a regular map, also denoted by~$f_i$, of the affine modular curve ${Y(N)=\Gamma(N)\backslash \poincare}$ to~$\C$, and our~$Z$ is the  image of  the map ${(f_1,\ldots, f_n):Y(N)\to\C^n}$. 

Furthermore, each~$f_i$ extends to a regular map ${\bar f_i:X(N)\to \P^1(\C)}$ of projective curves, where $X(N)$ is the standard compactification of $Y(N)$, as explained, for instance, in \cite[Section 2.4]{DS05}.  The image $\bar Z$ of the map ${(\bar f_1,\ldots, \bar f_n):X(N)\to\P^1(\C)^n}$ is Zariski closed in $\P^1(\C)^n$ and irreducible (being the image of an irreducible projective curve under a regular map). But for ${x\in X(N)}$ we have ${\bar f_i(x)=\infty}$ if and only if ${x\in X(N)\smallsetminus Y(N)}$ (we write ${\P^1(\C)=\C\cup\{\infty\}}$ in the obvious sense).  Hence ${Z=\bar Z\cap \C^n}$, which shows that~$Z$ is Zariski-closed in $\C^n$ and irreducible. This proves Part~\ref{izarclo}. 

Part~\ref{ispec2} is an easy consequence of Part~\ref{izarclo}. If~$f$ and~$g$ are two non-constant $j$-maps, then there exists~$N$ such that ${\Phi_N(f,g)=0}$. It follows that, up to coordinate permutations, any $j$-set is contained  in ${L_0\times M}$, where $L_0$ is a point whose coordinates are singular moduli and ${M\in \C^m}$ is defined as in~\eqref{eqm}. Since our $j$-set is irreducible and Zariski closed, it must be an irreducible component of ${L_0\times M}$, that is, a special curve. In particular, a $j$-set is an irreducible $1$-dimensional algebraic set defined over~$\bar\Q$. 

Conversely,  every special curve has (up to coordinate permutations) the shape ${L_0 \times \widetilde M}$, where~$\widetilde M$ is an irreducible component of a set ${M\subset \C^m}$ defined as in~\eqref{eqm}. Recall that two complex numbers $x,y$ satisfy ${\Phi_N(x,y)=0}$ if and only if~$x$ and~$y$ are $j$-invariants of two elliptic curves linked by a cyclic $N$-isogeny.  Now let ${(\xi_1,\ldots,\xi_m)}$ be a transcendental point\footnote{``transcendental'' means here that the coordinates of this point are not all algebraic over~$\Q$} of~$\widetilde M$. Then the numbers ${\xi_1, \ldots, \xi_m}$ are $j$-invariants of isogenous elliptic curves. Hence, if we write ${\xi_1=j(z)}$ with some ${z\in \poincare}$, then there exist ${\gamma_2,\ldots, \gamma_m\in \GL^+_2(\Q)}$ such that ${\xi_i=j(\gamma_iz)}$ for ${i=2,\ldots,m}$. 

Thus,~$\widetilde M$ shares a transcendental point with the $j$-set defined by the $j$-maps ${j, j\circ\gamma_2, \ldots, j\circ\gamma_m}$. Since both are Zariski-closed  irreducible $1$-dimensional algebraic sets defined over~$\bar\Q$, they must coincide. \qed

\bigskip

A similar ``parametric'' description can be given for higher dimensional special varieties. We do not go into this because we will not need it. 

\bigskip

Pila~\cite{Pi11} generalized the theorem of André by proving the following.

\begin{theorem}[Pila]
\label{thpila}
An algebraic set in $\C^n$ contains at most finitely many maximal special subvarieties. 
\end{theorem}

``Maximal'' is understood here in the set-theoretic sense: let~$V$ be an algebraic set in~$\C^n$ and ${M\subseteq V}$  a special variety;  we call~$M$ a \textsl{maximal special subvariety} of~$V$ if for any special variety~$M'$ such that ${M\subseteq M'\subseteq V}$ we have ${M=M'}$.

%(By an algebraic variety we mean an irreducible algebraic subset.)

If an algebraic curve is not special, than its only special subvarieties  are special points, and we recover the theorem of André.

\section{Main Lemma and Proof of Theorem~\ref{thcoli}}
\label{sml}

Theorem~\ref{thcoli} is an easy consequence of Pila's Theorem~\ref{thpila}  and the following lemma.

\begin{lemma}[``Main Lemma'']
\label{lmain}
Let $f_1,f_2,f_3,g_1,g_2,g_3$ be $j$-maps, not all constant. Assume that the determinant
{%\small
\begin{equation}
\label{ematfffggg}
\det\begin{bmatrix}
1&1&1\\
f_1&f_2&f_3\\
g_1&g_2&g_3
\end{bmatrix}
\end{equation}}%
is identically~$0$. Then at least one of the following holds: 
\begin{itemize}
\item
${f_1=f_2=f_3}$;
\item
${g_1=g_2=g_3}$;
\item
for some distinct ${k,\ell\in \{1,2,3\}}$ we have ${f_k=f_\ell}$ and ${g_k=g_\ell}$;
\item
${f_k=g_k}$ for  ${k=1,2,3}$.

\end{itemize} 
\end{lemma}

In this section we prove Theorem~\ref{thcoli} assuming the validity of the Main Lemma.
Lemma~\ref{lmain} itself will be proved in the subsequent sections.  

Consider the algebraic set in $\C^6$ consisting of the points ${(x_1,x_2,x_3,y_1,y_2,y_3)}$ satisfying
{%\small
\begin{equation}
\label{ev}
\begin{vmatrix}
1&1&1\\
x_1&x_2&x_3\\
y_1&y_2&y_3
\end{vmatrix}=0.
\end{equation}}%
Then Lemma~\ref{lmain} has the following consequence. 

\begin{corollary}
The algebraic set~\eqref{ev} has exactly six  maximal special subvarieties  of positive dimension: 
\begin{itemize}
\item
the subvariety $R_x$, defined in $\C^6$ by ${x_1=x_2=x_3}$;
\item
the subvariety $R_y$, defined in $\C^6$ by ${y_1=y_2=y_3}$;
\item
the three subvarieties $S_{k,\ell}$, defined in $\C^6$ by  ${x_k=x_\ell}$ and ${y_k=y_\ell}$, where ${k,\ell\in \{1,2,3\}}$ are distinct;
\item
the subvariety~$T$, defined in $\C^6$ by ${x_k=y_k}$  for  ${k=1,2,3}$.

\end{itemize} 
\end{corollary}

\paragraph{Proof}
Let~$\widetilde K$ be a special variety in (\ref{ev}) of positive dimension. By Proposition~\ref{pspec1} it contains a Zariski dense union of special curves. By Proposition~\ref{pspec2} each such curve is a $j$-set. By the Main Lemma each  $j$-set is contained  in one of the subvarieties above. The latter are clearly irreducible and also special; for example with $R_x$ we have ${n=6}$, ${d=4}$ and the partition with 
$$
l_0=0,\quad m_1=3,\quad m_2=m_3=m_4=1.
$$ 
Taking closures we see that $\widetilde K$ itself is also contained in one of them.
 \qed

\bigskip

Now we are ready to prove Theorem~\ref{thcoli}. Let 
$$
P_k=(x_k,y_k) \qquad (k=1,2,3)
$$
be three special points forming a collinear triple. Then the point ${Q=(x_1,x_2,x_3,y_1,y_2,y_3)}$ belongs to the algebraic set~\eqref{ev}. Moreover, since our points are pairwise distinct,~$Q$ does not belong to any of~$S_{k,\ell}$, and 
since the straight line passing through our points is not special,~$Q$ does not belong to any of $R_x,R_y,T$.

This shows that $\{Q\}$ is a zero-dimensional maximal special subvariety of the algebraic set~\eqref{ev}, and we complete the proof by applying Theorem~\ref{thpila}. \qed

\bigskip

The Main Lemma will be proved in Sections~\ref{sini}--\ref{smggngl}, after some preparations made in Sections~\ref{sroots}--\ref{sjmaps}.

\section{Roots of Unity}
\label{sroots}
In this section we collect some facts about roots of unity used in the proof of the Main Lemma.

\begin{lemma}
\label{lsumroots}
Let~$\alpha$ be a sum of~$k$ roots of unity, and~$N$ a non-zero integer. Assume that ${N\mid \alpha}$ (in the ring of algebraic integers). Then either ${\alpha=0}$ or ${k\ge |N|}$. 
\end{lemma}

\paragraph{Proof}
Assume ${\alpha\ne 0}$ and write ${\alpha=N\beta}$, where~$\beta$ is a non-zero algebraic integer. Then there exists an embedding ${\Q(\alpha)\stackrel\sigma\to\C}$ such that ${|\beta^\sigma|\ge 1}$. It follows that ${|N|\le |\alpha^\sigma|}$. But, since~$\alpha$ is a sum of~$k$ roots of unity, we have ${|\alpha^\sigma|\le k}$. \qed 

\begin{lemma}
\label{lquadratic}
Let $a,b$ be non-zero rational numbers and $\eta,\theta$ roots of unity. Assume that ${\alpha=a\eta+b\theta}$ is of degree~$1$ or~$2$ over~$\Q$. Then $\Q(\alpha)$ is one of the fields~$\Q$, $\Q(i)$, $\Q(\sqrt{-2})$, $\Q(\sqrt{-3})$,  $\Q(\sqrt2)$, $\Q(\sqrt3)$, $\Q(\sqrt5)$, and after a possible swapping of $a\eta$ and $b\theta$, and possible replacing of $(a,\eta)$ by $(-a,-\eta)$ and/or $(b,\theta)$ by $(-b,-\theta)$, we have the following.

\begin{enumerate}
\item
\label{iq}
If ${\Q(\alpha)=\Q}$ then:

\begin{enumerate}
\item
\label{iipmone}

either both~$\eta$ and~$\theta$ are $\pm1$,

\item
\label{iicubic}
or~$\eta$ is a primitive cubic root of unity, ${\theta=\eta^{-1}}$ and ${a=b}$,

\item
\label{iizero}
or ${\theta=-\eta}$  and ${a=b}$.

\end{enumerate}

\item
If ${\Q(\alpha)=\Q(i)}$ then:

\begin{enumerate}

\item
\label{iiione}

either ${\eta=i}$  and ${\theta\in \{1,i\}}$, %one of ${\eta,\theta}$ is $\pm i$ and the other is $\pm i$ or $\pm1$,

\item
\label{i12i}
or~$\eta$ is a primitive $12$th root of unity, ${\theta = -\eta^{-1}}$ and ${a=b}$. %${b\theta = -a\eta^{-1}}$.

\end{enumerate}

\item
\label{izetathree}
If ${\Q(\alpha)=\Q(\sqrt{-3})}$ then~$\eta$   %one of ${\eta,\theta}$ is $\pm\zeta_3$, where~$\zeta_3$ 
is a primitive cubic root of unity, and~$\theta$ is a cubic root of unity (primitive or not). %the other is $\pm \zeta_3$ or $\pm1$.

\item
\label{i8-2}
If ${\Q(\alpha)=\Q(\sqrt{-2})}$ then~$\eta$ is a primitive $8$th  root of unity, ${\theta = -\eta^{-1}}$ and ${a=b}$. % and ${b\theta = -a\eta^{-1}}$. 

\item
\label{i82}
If ${\Q(\alpha)=\Q(\sqrt{2})}$ then~$\eta$ is a primitive $8$th root of unity, ${\theta = \eta^{-1}}$ and ${a=b}$.  %and ${b\theta = a\eta^{-1}}$,

\item

If ${\Q(\alpha)=\Q(\sqrt{3})}$ then~$\eta$ is a primitive $12$th root of unity, and

\begin{enumerate}
\item
\label{i121}
either ${\theta = \eta^{-1}}$,\quad ${a=b}$, %  and ${b\theta = a\eta^{-1}}$,  

\item
\label{i122}
or ${\theta=-\eta^3(=\pm i)}$,\quad ${a=2b}$.   
%, after a possible swapping of $a\eta$ and $b\theta$, we have the following, and ${2b\theta=-a\eta^3}$; in particular, ${b=\pm a/2}$ and ${\theta=\pm i}$. 

\end{enumerate}

\item
\label{izetafive}
If ${\Q(\alpha)=\Q(\sqrt{5})}$ then~$\eta$ is a primitive $5$th root of unity,  ${\theta = \eta^{-1}}$ and ${a=b}$. % ${b\theta = a\eta^{-1}}$. 

\end{enumerate}
\end{lemma}

\paragraph{Proof}
Without loss of generality we may assume that~$a$ and~$b$ are coprime integers. Let~$N$ be the order of the multiplicative group generated by~$\eta$ and~$\theta$, and ${L=\Q(\eta,\theta)}$; then 
${[L:\Q]=\ph(N)}$, where~$\ph$ is  Euler's totient function.  

If ${\ph(N)\le2}$ then ${N\in \{1,2,3,4,6\}}$, and we have one of the options~\ref{iq},~\ref{iiione} or~\ref{izetathree}. 
If ${\alpha=0}$ then we have option~\ref{iizero}.

From now on we assume that ${\ph(N)>2}$ and ${\alpha\ne 0}$. Since ${\ph(N)>2}$, there exists  ${\sigma\in \Gal(L/\Q)}$ such that ${(\eta^\sigma,\theta^\sigma)\ne(\eta,\theta)}$, but ${\alpha^\sigma=\alpha}$. We obtain 
\begin{equation}
\label{ebothsides}
a(\eta-\eta^\sigma)=b(\theta^\sigma-\theta). 
\end{equation}
By our choice of~$\sigma$ both sides of~\eqref{ebothsides} are non-zero. Since~$a$ and~$b$ are coprime integers, we have ${a\mid(\theta^\sigma-\theta)}$, whence ${|a|\le 2}$ by Lemma \ref{lsumroots}. Similarly, 
${|b|\le 2}$. It follows that ${(a,b)\in \{(\pm1,\pm1),(\pm1,\pm2),(\pm2,\pm1)\}}$. 
Swapping (if necessary) $a\eta$ and $b\theta$, and  replacing (if necessary) $(a,\eta)$ by $(-a,-\eta)$ and/or $(b,\theta)$ by $(-b,-\theta)$,  we may assume that ${a\in \{1,2\}}$ and  ${b=1}$. The rest of the proof splits into two cases. 

\bigskip

\noindent
\underline{The case ${a=2}$, ${b=1}$}\quad
In this  case~\eqref{ebothsides} writes as ${2(\eta-\eta^\sigma)=\theta^\sigma-\theta}$.  We must have ${\theta^\sigma=-\theta}$; otherwise all the conjugates of the non-zero algebraic integer ${(\theta^\sigma-\theta)/2}$ would be of absolute value strictly smaller than~$1$. Thus, we obtain ${\eta-\eta^\sigma+\theta=0}$. Three roots of unity may sum up to~$0$ only if they are proportional to $(1,\zeta_3,\zeta_3^{-1})$, where~$\zeta_3$ is a primitive cubic of unity.   We obtain ${\theta/\eta=\zeta_3^{-1}}$, and ${\eta=\alpha(a+b\zeta_3^{-1})^{-1}}$ is of degree at most~$4$ over~$\Q$. Since ${\theta=\eta^\sigma-\eta\in \Q(\eta)}$, we obtain ${L=\Q(\eta)}$; in particular,~$\eta$ is a primitive $N$th root of unity.

Thus, ${\ph(N)=[\Q(\eta):\Q]\le 4}$, and in fact ${\ph(N)=4}$ because  ${\ph(N)>2}$. Since ${-\eta^\sigma/\eta=\zeta_3}$, we must have ${3\mid N}$. Together with ${\ph(N)=4}$ this implies that  ${N=12}$ and~$\eta$ is a primitive $12$th root of unity. Hence we have the option~\ref{i122}.  

\bigskip

\noindent
\underline{The case ${a=b=1}$} \quad   In this case ${\eta-\eta^\sigma+\theta-\theta^\sigma=0}$. Four roots of unity may sum up to~$0$ only if two of them sum up to~$0$ (and the other two sum up to~$0$ as well). Since ${\eta\ne \eta^\sigma}$ and ${\eta\ne-\theta}$ (because ${\alpha\ne 0}$), we have  ${\eta=\theta^\sigma}$ and ${\eta^\sigma=\theta}$. This implies that ${L=\Q(\eta)=\Q(\theta)}$, both~$\eta$ and~$\theta$ are primitive $N$th roots of unity,  and ${\sigma^2=1}$.

We claim that the subgroup ${H=\{1,\sigma\}}$ is the stabilizer of $\Q(\alpha)$ in ${G=\Gal(L/\Q)}$. 
Thus, let ${\varsigma\in G}$ satisfy ${\alpha^\varsigma=\alpha}$. 
Since ${\eta+\eta^\sigma-\eta^\varsigma-\eta^{\sigma\varsigma}=0}$ and ${\eta+\eta^\sigma\ne0}$, we must have either ${\eta=\eta^\varsigma}$ or ${\eta=\eta^{\sigma\varsigma}}$. Since ${L=\Q(\eta)}$, in the first case we have ${\varsigma=1}$ and in the second case ${\varsigma=\sigma^{-1}=\sigma}$. 

Thus,~$H$ is the stabilizer of $\Q(\alpha)$. Since ${|H|=2}$ and ${[G:H]=[\Q(\alpha):\Q]=2}$, we obtain ${\ph(N)=|G|=4}$, which implies that ${N\in\{5,8,10,12\}}$.

Now if ${N=5}$ then we have option~\ref{izetafive}. If ${N=10}$ then, replacing  $(a,\eta)$ by $(-a,-\eta)$ and $(b,\theta)$ by $(-b,-\theta)$, we obtain option~\ref{izetafive} as well. If ${N=8}$ then we have one of the options~\ref{i8-2} or~\ref{i82}. Finally, if ${N=12}$ then we have one of the options~\ref{i12i} or~\ref{i121}. \qed

\section{Singular Moduli}
\label{smodi}

In this section  we collect miscellaneous  properties of singular moduli used in the sequel. We start by recalling 
 the notion of the \textsl{discriminant} of a singular modulus. Let ${\tau\in \poincare}$ be algebraic of degree~$2$; the endomorphism ring of the lattice $\Z\tau+\Z$ is an order in the imaginary quadratic field $\Q(\tau)$; the discriminant ${\Delta=\Delta_\tau}$ of this order will be called the \textsl{discriminant} of the singular modulus $j(\tau)$. This discriminant is a negative integer satisfying ${\Delta\equiv 0,1\bmod 4}$. %; conversely, any integer ${\Delta<0}$ satisfying ${\Delta\equiv 0,1\bmod 4}$ is a discriminant of some singular modulus; for instance, of the singular modulus ${j\bigl((\Delta+\sqrt\Delta)/2\bigr)}$; this latter will be called the \textsl{principal} singular modulus of discriminant~$\Delta$. 

It is well-known (see, for instance, \cite[Section~11]{Co89}) that 
\begin{itemize}
\item
any singular modulus of discriminant~$\Delta$ is an algebraic integer of degree equal to the class number of~$\Delta$, denoted $h(\Delta)$;

\item
the  singular moduli of discriminant~$\Delta$ are all conjugate over~$\Q$; moreover, they form a complete set of $\Q$-conjugates.  

\end{itemize}

A full description of singular moduli of given discriminant~$\Delta$ is well-known as well. Denote by ${T=T_\Delta}$ the set of triples of integers $(a,b,c)$ such that 
\begin{equation*}
%\label{ekuh}
\begin{gathered}
\gcd(a,b,c)=1, \quad \Delta=b^2-4ac,\\
\text{either\quad $-a < b \le a < c$\quad or\quad $0 \le b \le a = c$}
\end{gathered}
\end{equation*}
\label{pallconj}
Then the map
\begin{equation}
\label{ejabcmap}
(a,b,c)\mapsto j\left(\frac{b+\sqrt{\Delta}}{2a}\right)
\end{equation}
defines a bijection from $T_\Delta$ onto the set of singular moduli of discriminant~$\Delta$. 
In particular, ${h(\Delta)=|T_\Delta|}$.
The proof of this  is a compilation of several classical facts, some of which go back to Gauss; see, for instance,   \cite[Section~2.2]{BLP16} and the references therein. 

It is crucial for us that the set~$T_\Delta$ has only one triple $(a,b,c)$ with ${a=1}$. The corresponding singular modulus will be called the \textsl{principal} singular modulus of discriminant~$\Delta$. 
Note that the principal singular modulus is a real number; in particular,
\begin{equation}
\label{ereconj}
\text{any singular modulus has a real $\Q$-conjugate}. 
\end{equation}

There exist exactly 13 discriminants~$\Delta$ with ${h(\Delta)=1}$. The corresponding singular moduli (and only them) are rational integers. The full list of the 13 rational singular moduli is well-known and reproduced in Table~\ref{talist1}. 
\begin{table}
\caption{Discriminants~$\Delta$ with ${h(\Delta)=1}$ and the corresponding singular moduli}
\label{talist1}
{\footnotesize
\begin{equation*}
%\label{elist1}
\begin{array}l
\begin{array}{r|lllllllll}
\Delta
&-3&-4&-7&-8&-11&-12&-16&-19&-27\\
%j(\tau_\Delta)
j&0&1728&-3375&8000& -32768& 54000& 287496& -884736& -12288000\\
\end{array}\\
\hline
\begin{array}{r|lllllllll}
\Delta&-28&-43&-67&-163\\
%j(\tau_\Delta)
j&16581375&-884736000&
-147197952000& -262537412640768000
\end{array}
\end{array}
\end{equation*}  
}
\end{table}

Finally, we use the inequality
\begin{equation}
\label{etwothousand}
\bigl||j(\tau)|-e^{2\pi\Im\tau}\bigr|\le 2079 
\end{equation}
which holds for every ${\tau\in \poincare}$ satisfying ${\Im\tau\ge \sqrt3/2}$ (see \cite[Lemma~1]{BMZ13}). In particular, if ${(a,b,c)\in T_\Delta}$ then the number 
$$
\tau(a,b,c)=\frac{b+\sqrt\Delta}{2a}
$$
satisfies ${\Im\tau(a,b,c)\ge \sqrt3/2}$ (see \cite[page~403, equation~(8)]{BLP16}). Hence~\eqref{etwothousand} applies with ${\tau=\tau(a,b,c)}$. 

All the facts listed above will be repeatedly used in this section, sometimes without a special reference.

\begin{lemma}
\label{lverybig}
Let~$x$ be a  singular modulus and let~$x'$ be the principal singular modulus of the same discriminant. Then either ${x=x'}$ or ${|x'|>|x|+180000}$. 
\end{lemma}

\paragraph{Proof}
Let~$\Delta$ be the common discriminant of~$x$ and~$x'$. We may assume that ${|\Delta|\ge 15}$, otherwise ${h(\Delta)=1}$ and there is nothing to prove. We assume that ${x\ne x'}$ and will use~\eqref{etwothousand} to estimate $|x|$ from above and $|x'|$ from below. 

We have ${x=j(\tau)}$ and ${x'=j(\tau')}$, where ${\tau =\tau(a,b,c)}$ and ${\tau'=\tau(a',b',c')}$ for some ${(a,b,c),(a',b',c')\in T_\Delta}$. Since~$x'$ is principal, and~$x$ is not, we have ${a'=1}$ and ${a\ge 2}$. Hence 
$$
\Im\tau'=\pi|\Delta|^{1/2}, \qquad \Im\tau =\frac{\pi|\Delta|^{1/2}}a\le \frac{\pi|\Delta|^{1/2}}2.
$$
We obtain 
$$
|x'|\ge e^{\pi|\Delta|^{1/2}}-2079, \qquad |x|\le e^{\pi|\Delta|^{1/2}/2}+2079, 
$$
which implies
$$
|x'|-|x|\ge e^{\pi|\Delta|^{1/2}}-e^{\pi|\Delta|^{1/2}/2}-4158\ge e^{\pi\sqrt{15}}-e^{\pi\sqrt{15}/2}-4158 > 180000,
$$
as wanted. \qed

\begin{lemma}
\label{lrootofone}
%The sum of two singular moduli cannot be equal to~$0$ or $1488$. 
Let~$x$,~$y$ be singular moduli and let ${a,b\in \Z}$ be such that ${|a|,|b|\le 90000}$. Assume that ${y\ne b}$ and that ${(x-a)/(y-b)}$ is a root of unity. Then either ${x=y}$ or ${x,y\in \Z}$.
In particular, if $x/y$ is a root of unity (with ${y\ne 0}$) or if ${(x-744)/(y-744)}$ is a root of unity then ${x=y}$. 
\end{lemma}

\paragraph{Proof}
Let~$x'$ and~$y'$ be the principal singular moduli of the same discriminants as~$x$ and~$y$. We may assume that ${|x'|\ge |y'|}$.  We may further assume, by conjugating, that ${x=x'}$. Then ${y=y'}$  as well, since otherwise ${|y|<|y'|-180000}$ by Lemma~\ref{lverybig}, and we obtain 
$$
|y|+90000\ge |y-b|=|x-a|=|x'-a|\ge|x'|-90000 \ge |y'|-90000>|y|+90000,
$$ 
a contradiction. Thus, both~$x$ and~$y$ are principal singular moduli. In particular,  both are real, which implies ${x-a=\pm (y-b)}$. 

Now Theorem~1.2 of~\cite{ABP15} implies one of the following options:
\begin{enumerate}
\item
${x=y}$ and ${a=b}$;
 
\item
${x,y\in\Z}$;

\item
\label{idegtwo}
$x$ and~$y$ are distinct and of degree~$2$ over~$\Q$.
\end{enumerate}

We have to rule out option~\ref{idegtwo}.   Thus, assume it to be the case and let  ${f(T)=T^2+AT+C}$ and ${g(T)=T^2+BT+D}$ be the $\Q$-minimal polynomials of~$x$ and~$y$. Since~$x$ and~$y$ are both principal and distinct, they are not $\Q$-conjugate, which means that the polynomials~$F$ and~$G$ are distinct. We have either ${x+y=a+b}$ or ${x-y=a-b}$. Taking $\Q$-traces, we obtain ${A+B=2(a+b)}$ or ${A-B=2(a-b)}$. In particular, we have either ${|A+B|\le 360000}$ or ${|A-B|\le 360000}$.

However, our~$F$ and~$G$ are among the~29 Hilbert class polynomials associated to the imaginary quadratic orders of class number~$2$. The full list of such polynomials can be found in Table~2 of~\cite{BLP16}. 
A quick inspection of this table shows that if~$A$ and~$B$ are middle coefficients of two distinct polynomials from this table, then ${|A+B|>360000}$ and ${|A-B|>360000}$. Hence option~\ref{idegtwo} is impossible. This proves the first statement of the lemma.

In the special cases ${a=b=0}$ or ${a=b=744}$ we must have either ${x=y}$ or 
\begin{equation}
\label{eimposs}
x,y\in \Z, \quad x\ne y , \quad x+y\in \{0,1488\}. 
\end{equation}
Inspecting Table~\ref{talist1}, we find out that~\eqref{eimposs} is impossible.  The lemma is proved. \qed

\begin{lemma}
\label{lmonoton}
Let~$x$ and~$y$ be distinct principal singular moduli. Then ${\bigl||x|-|y|\bigr|>1600}$. 
\end{lemma}

\paragraph{Proof}
Denote by~$\Delta_x$ and~$\Delta_y$ the discriminants of~$x$ and~$y$, respectively. We will assume that ${|\Delta_x|>|\Delta_y|}$. If ${|\Delta_x|\le 12}$ then ${h(\Delta_x)=1}$, and the statement follows by inspection of Table~\ref{talist1}. And if ${|\Delta_x|\ge15}$ then 
\begin{align*}
|x|-|y|&\ge (e^{\pi|\Delta_x|^{1/2}}-2079)-(e^{\pi|\Delta_y|^{1/2}}+2079)\\
&\ge e^{\pi|\Delta_x|^{1/2}}- e^{\pi|\Delta_x-1|^{1/2}}-4158\\ 
&\ge e^{\pi\sqrt{15}}-e^{\pi\sqrt{14}}-4158\\
&>60000,
\end{align*}
which is much stronger than needed. The lemma is proved. \qed

\begin{lemma}
\label{ltwoelem}
Let~$x$ be a singular modulus, and assume that the number field $\Q(x)$ is a Galois extension of~$\Q$. Then the Galois group of ${\Q(x)/\Q}$ is $2$-elementary; that is, isomorphic to 
${(\Z/2\Z)^k}$ for some~$k$. 
\end{lemma}

\paragraph{Proof}
This is well-known; see, for instance, Corollary~3.3 from~\cite{ABP15}. \qed

\begin{lemma}
\label{lcute}
Let $x,y$ be singular moduli and ${\eps,\eta}$ roots of unity. Then ${\eps (x-744)+\eta (y-744)}$ is not a root of unity.  
\end{lemma}

\paragraph{Proof}
We will assume that 
$$
\eps (x-744)+\eta (y-744)=1
$$
and derive a contradiction. We have clearly
\begin{equation}
\label{ecolumbus}
\bigl||y|-|x|\bigr|\le 1489.
\end{equation}
We follow the same strategy as in the proof of Lemma~\ref{lrootofone}. We denote by~$x'$ and~$y'$ the principal moduli of the same discriminants as~$x$ and~$y$, respectively, and we may assume that ${|x'|\ge|y'|}$ and ${x=x'}$. We claim that ${y=y'}$ as well. Indeed, if ${y\ne y'}$ then Lemma~\ref{lverybig} implies that 
$$
|y|+1489\ge |x| =|x'|\ge |y'|> |y|+180000,
$$ 
a contradiction.

Thus, we may assume that both~$x$ and~$y$ are principal singular moduli. Lemma~\ref{lmonoton} and inequality~\eqref{ecolumbus} imply that ${x=y}$. Thus,
$$
(\eps+\eta) (x-744)=1.
$$
In particular ${0\ne\eps+\eta\in \R}$, which implies ${\eta=\eps^{-1}}$. 

Lemma~\ref{ltwoelem} implies that the Galois group of the number field ${\Q(x)=\Q(\eps+\eps^{-1})}$ is $2$-elementary. Since ${\Q(\eps+\eps^{-1})}$ is a subfield of degree at most~$2$ in $\Q(\eps)$, the Galois group of ${\Q(\eps)/\Q}$ is either $2$-elementary or $\Z/4\Z$ times a $2$-elementary group. But this group is $(\Z/n\Z)^\times$, where~$n$ is the order of the root of unity~$\eps$. Using the well-known structure of the multiplicative group $(\Z/n\Z)^\times$ (see, for instance, \cite[Theorem~3 in Section 4.1]{IR90}), one easily finds out that any integer~$n$ with the property ``the group $(\Z/n\Z)^\times$ is either $2$-elementary or $\Z/4\Z$ times a $2$-elementary group'' divides either $48$ or $120$.  It follows that 
${|\eps+\eps^{-1}|\ge 2\sin(\pi/60)}$ (recall that ${\eps+\eps^{-1}=\eps+\eta\ne 0}$).
Hence 
$$
|x-744|\le \frac{1}{2\sin(\pi/60)}< 10. 
$$
No principal singular modulus satisfies the latter inequality. The lemma is proved. \qed 

%For the proof see \cite[Lemma~3.5]{ABP15}. In \cite{ABP15} it is assumed that the discriminant is at least~$11$; but for ${\Delta<11}$ we have ${h(\Delta)=1}$, so there is only one singular modulus, which is principal.

%In particular, we have the following. 

\begin{lemma}
\label{ljne744etc}
The numbers  
${744,\ 744\pm1,\ 744\pm2,\ 
744\pm196884,\ 744\pm1\pm196884,}  \ 744\pm2\cdot196884$
are not singular moduli. 
\end{lemma}

\paragraph{Proof} The proof is just by inspection of Table~\eqref{talist1}. \qed

\begin{lemma}
\label{ljneroot}
Let~$\theta$ be a root of unity. Then ${744+\theta}$ and ${744+196884\theta}$ are not  singular moduli.
\end{lemma}

\paragraph{Proof}
%This is an immediate consequence of Lemmas~\ref{ljnecyclo} and~\ref{ljne744etc}.\qed
If  ${744+\theta}$ or ${744+196884\theta}$ is a singular modulus, then the cyclotomic field $\Q(\theta)$ has a real embedding by~\eqref{ereconj}, which is possible only if ${\theta=\pm1}$. Now apply Lemma~\ref{ljne744etc}.\qed

\begin{lemma}
\label{lpluricyc}
Assume that a singular modulus of discriminant~$\Delta$ is a sum of~$k$ roots of unity. Then 
$$
|\Delta|\le \pi^{-2}(\log(k+2079))^2. 
$$
\end{lemma}

\paragraph{Proof}
We may assume that our modulus (denote it~$x$) is principal, and, 
as in the proof of Lemma~\ref{lverybig}, deduce from this that it satisfies ${|x|\ge e^{\pi|\Delta|^{1/2}}-2079}$. On the other hand, since~$x$ is a sum of~$k$ roots of unity, we have ${|x|\le k}$. Whence the result. \qed

%Replacing our singular modulus by a conjugate, we may assume that it is the principal modulus $j(\tau)$ with ${\tau=(\Delta+\sqrt\Delta)/2}$. Using the inequality ${||j(\tau)|-q^{-1}|\le 2079}$, where, as usual, ${q=e^{2\pi i \tau}}$ (see \cite[Lemma~1]{BMZ13}), we find ${|j(\tau)|\ge e^{\pi|\Delta|^{1/2}}-2079}$. On the other hand, since $j(\tau)$ is a sum of~$k$ roots of unity, we have ${|j(\tau)|\le k}$. Whence the result. \qed

\begin{lemma}
\label{ljnetworoots}
Let ${\eta,\theta}$ be roots of unity,~$x$ a singular modulus and ${a,b,c\in \Z}$.  % of discriminant~$\Delta$. 
Assume that 
$$
x=a\eta+b\theta+c,  \quad a,b\ne 0, \quad |a|+|b|+|c|\le 3400000. 
$$
Then one of the following options takes place.

\begin{itemize}
\item
We have ${x\in \Z}$.

\item
After possible replacing of $(a,\eta)$ by $(-a,-\eta)$ and/or $(b,\theta)$ by $(-b,-\theta)$, we have the following:~$\eta$ is a primitive $5$th root of unity,   ${\theta = \eta^{-1}}$,  ${a=b}$ and 
%${(a,b,c)}$ is one of the following
\begin{equation}
\label{eacs}
(a,c) \in \{(85995,-52515), (-85995,-138510),
(565760,914880), (-565760,349120)\}.
\end{equation}
\end{itemize}

\end{lemma}%${\Delta=-15}$,\quad ${j(\tau)=(-191025\pm85995\sqrt5)/2}$,\quad 

\paragraph{Proof}
Let~$\Delta$ be the discriminant of the singular modulus~$x$. Lemma~\ref{lpluricyc} implies that 
\begin{equation}
\label{edelta<}
|\Delta|\le \pi^{-2}(\log(3400000+2079))^2<22.92.
\end{equation}
Assume that ${x\notin\Z}$; then ${h(\Delta)>1}$. Among negative quadratic discriminants satisfying~\eqref{edelta<} all but two have class number~$1$; these two are    ${\Delta=-15}$ and ${\Delta=-20}$. In both cases ${h(\Delta)=2}$ and ${\Q(x)=\Q(\sqrt5)}$, so option~\ref{izetafive} of Lemma~\ref{lquadratic} applies in both cases. After possible replacing of $(a,\eta)$ by $(-a,-\eta)$ and/or $(b,\theta)$ by $(-b,-\theta)$, we obtain the following:~$\eta$ is a primitive $5$th root of unity,  ${\theta = \eta^{-1}}$ and ${a=b}$, so we have 
${x=a(\eta+\eta^{-1})+c}$.

The two singular moduli of discriminant ${\Delta=-15}$ are  
{%\small
$$
\frac{-191025\pm85995\sqrt5}2= -\frac{191025}2\pm85995\left(\frac12+\eta+\eta^{-1}\right)=
\begin{cases}
\text{either}& \hphantom{-}85995(\eta+\eta^{-1})-52515, \\ 
\text{or}& -85995(\eta+\eta^{-1})-138510,
\end{cases}
$$}%
which gives us the first two options in~\eqref{eacs}

Similarly, the two singular moduli of discriminant ${\Delta=-20}$ are   ${632000\pm282880\sqrt5}$, which gives the other two options. \qed

\section{Rational Matrices}
\label{srmat}

In this section we obtain some elementary properties of $\Q$-matrices, which will be used in our study of $j$-maps in Section~\ref{sjmaps}.

Recall that we denote by $\GL_2^+(\Q)$ the subgroup of $\GL_2(\Q)$ consisting of matrices of positive determinant. 
Unless the contrary is stated explicitly, in this section \textsl{matrix} refers to an element in $\GL_2^+(\Q)$.  We call two matrices~$A$ and~$A'$ \textsl{equivalent} (notation: ${A\sim A'}$) if there exists a matrix ${B\in \SL_2(\Z)}$  and a scalar ${\lambda\in \Q^\times}$ such that ${A'=\lambda BA}$. 

For ${a,b\in \Q}$ we define $\gcd(a,b)$ as the non-negative ${\delta\in \Q}$ such that ${a\Z+b\Z=\delta\Z}$. 

Given a matrix ${A=\begin{bsmallmatrix}
a&b\\c&d
\end{bsmallmatrix}}$, we define the \textsl{normalized left content} of~$A$ by 
$$
\nlc(A)= \frac{\gcd(a,c)^2}{\det A}. 
$$
Clearly, ${\nlc(A)=\nlc(A')}$ if ${A\sim A'}$.

\begin{proposition}
\label{ptri}
Every matrix~$A$ is equivalent to an upper triangular matrix of the form ${\begin{bsmallmatrix}
a&b\\
0&1
\end{bsmallmatrix}}$ with ${a>0}$, where ${a=\nlc(A)}$.   We have ${\begin{bsmallmatrix}
a&b\\
0&1
\end{bsmallmatrix}\sim \begin{bsmallmatrix}
a'&b'\\
0&1
\end{bsmallmatrix}}$ if and only if ${a=a'}$ and ${b\equiv b'\bmod \Z}$. 
\end{proposition}

\paragraph{Proof}
It suffices to show that~$A$ is equivalent to an upper triangular matrix;  the rest is easy. Let $\begin{psmallmatrix}
x\\y
\end{psmallmatrix}$ 
be the left column of~$A$ and ${\delta=\gcd(x,y)}$. Then ${x/\delta,y/\delta\in \Z}$ and there exist ${u,v\in \Z}$ such that ${ux+vy=\delta}$. Multiplying~$A$ on the left by the matrix ${\begin{bsmallmatrix*}[l]
\hphantom{-}u&v\\
-y/\delta&x/\delta
\end{bsmallmatrix*}\in \SL_2(\Z)}$, we obtain an upper triangular matrix.\qed

{\sloppy

\begin{proposition}
\label{ptwo}
Let $A_1,A_2$ be non-equivalent matrices. Then there exists a matrix~$B$  such that ${\nlc(A_1B)\ne \nlc(A_2B)}$. 
\end{proposition}

\paragraph{Proof}
We may assume that ${\nlc(A_1)=\nlc(A_2)}$ (otherwise there is nothing to prove). Multiplying  on the right by $A_1^{-1}$, we may assume that ${A_1=\begin{bsmallmatrix}
1&0\\
0&1
\end{bsmallmatrix}}$.  We may further assume that ${A_2=\begin{bsmallmatrix}
a&b\\
0&1
\end{bsmallmatrix}}$. Since ${a=\nlc(A_2)=\nlc(A_1)=1}$, we have %${a=1}$. Thus, 
${A_2=\begin{bsmallmatrix}
1&b\\
0&1
\end{bsmallmatrix}}$, where ${b\notin \Z}$ since ${A_2\nsim A_1}$.

}

Now ${B=\begin{bsmallmatrix}
1&0\\
-b^{-1}&1
\end{bsmallmatrix}}$ would do. Indeed, 
$$
\nlc(A_1B)=\nlc(B)= \gcd(-b^{-1},1)^2, \quad \nlc(A_2B)=\nlc\left(\begin{bsmallmatrix}
0&b\\
-b^{-1}&1
\end{bsmallmatrix}\right)= b^{-2}, 
$$
and we have to prove that ${\gcd(-b^{-1},1)\ne |b|^{-1}}$. This is equivalent to ${\gcd(1,b)\ne 1}$, which is true because ${b\notin\Z}$. \qed

\bigskip

One may wonder if the same statement holds true for more than two  matrices: \textsl{given pairwise non-equivalent matrices ${A_1,\ldots, A_n}$, does there exists a matrix ${B\in \GL_2^+(\Q)}$ such that 
${\nlc(A_1B), \ldots, \nlc(A_nB)}$
are pairwise distinct?}  The proof of the  Main Lemma could have been drastically simplified if it were the case. Unfortunately, the answer is ``no'' already for three matrices, as the following example shows.

\begin{example}
Let 
{%\small
$$
A_1= 
\begin{bmatrix}
1&0\\
0&1
\end{bmatrix}, \quad 
A_2= 
\begin{bmatrix}
1&1/2\\
0&1
\end{bmatrix}, \quad 
A_3= 
\begin{bmatrix}
4&0\\
0&1
\end{bmatrix}. 
$$}%
We claim that for any matrix~$B$, at least two of the numbers 
$$
\nlc(A_1B),\quad \nlc(A_2B),\quad \nlc(A_3B)
$$ 
are equal. Indeed, write ${B=\begin{bsmallmatrix}
a&b\\c&d
\end{bsmallmatrix}}$. After multiplying by a suitable scalar, we may assume that ${c=2}$. Now
{%\small
$$
\nlc(A_1B)= \frac{\gcd(a,2)^2}{\det B}, \quad
\nlc(A_2B)= \frac{\gcd(a+1,2)^2}{\det B}, \quad
\nlc(A_3B)= \frac{\gcd(4a,2)^2}{4\det B},
$$}%
and we must show that among the three numbers
$$
\gcd(a,2), \quad \gcd(a+1,2), \quad \frac12\gcd(4a,2),
$$
there are two  equal. And this is indeed the case:
\begin{itemize}
\item
if ${\ord_2(a)>0}$ then ${\gcd(a+1,2)= \frac12\gcd(4a,2)}$; 
\item
if ${\ord_2(a)=0}$ then ${\gcd(a,2)= \frac12\gcd(4a,2)}$;

\item
if ${\ord_2(a)<0}$ then ${\gcd(a,2)= \gcd(a+1,2)}$. \qed
\end{itemize} 
\end{example}

Still, it is possible to prove something. 

\begin{proposition}
\label{pthree}
Let $A_1,A_2,A_3$ be pairwise non-equivalent matrices. Then there exists a matrix~$B$  such that among the numbers $\nlc(A_1B)$,  $\nlc(A_2B)$, $\nlc(A_3B)$, one is strictly bigger than the two others. 
\end{proposition}

\paragraph{Proof}
We may assume that ${A_k=
\begin{bsmallmatrix}a_k&\ast\\0&1\end{bsmallmatrix}}$ for ${k=1,2,3}$. 
If the numbers $a_k$ are pairwise distinct then there is nothing to prove. Hence we may assume that ${a_1=a_2}$. Multiplying  on the right by $A_3^{-1}$, and afterwards by a suitable diagonal matrix,   we may assume that 
{%\small
$$
A_1=
\begin{bmatrix}
1&b_1\\
0&1
\end{bmatrix},\quad
A_2=
\begin{bmatrix}
1&b_2\\
0&1
\end{bmatrix},\quad
A_3=
\begin{bmatrix}
a^{-1}&0\\
0\hphantom{{}^{-1}}&1
\end{bmatrix},
$$}%
where ${a>0}$. %We may also assume that ${b_1,b_2>0}$, adding to each a suitable integer. 
Since ${A_1\nsim A_2}$ we have ${b_1\not\equiv b_2\bmod \Z}$, and we may assume that 
${b_1\notin\Z}$. 

Set ${B=\begin{bsmallmatrix}
1&0\\
-b_1^{-1}&1
\end{bsmallmatrix}}$. Then 
\begin{equation}
\label{ethreenum}
\nlc(A_1B)= b_1^{-2}, \quad
\nlc(A_2B)= \gcd(1-b_1^{-1}b_2, b_1^{-1})^2, \quad
\nlc(A_3B)= a\gcd(a^{-1},b_1^{-1})^2. 
\end{equation}
Multiplying numbers~\eqref{ethreenum} by ${ab_1^2}$, we must show that among the three numbers 
\begin{equation}
\label{ethreenumbis}
a, \quad
a \gcd(b_1-b_2, 1)^2, \quad
\gcd(b_1,a)^2. 
\end{equation}
one is strictly bigger than the others.

If the numbers in~\eqref{ethreenumbis} are pairwise distinct then there is nothing to prove. Now assume that two of them are equal. Since ${b_1\not\equiv b_2\bmod \Z}$, we have ${\gcd(b_1-b_2, 1)<1}$, and, in particular, the first two of them % numbers~\eqref{ethreenumbis} 
are distinct. 

Further, the equality ${a=\gcd(b_1,a)^2}$ is not possible either. Indeed, in this case for any prime number~$p$ we would have 
$$
\ord_p(a)=2\min\{\ord_p(a), \ord_p(b_1)\},
$$ 
which implies that either ${\ord_p(a)=2\ord_p(b_1)>0}$ or ${\ord_p(b_1)\ge \ord_p(a)=0}$. In particular, ${\ord_p(b_1)\ge 0}$ for any~$p$, contradicting our assumption ${b_1\notin \Z}$. 

Thus, the only possibility is ${a \gcd(b_1-b_2, 1)^2=
\gcd(b_1,a)^2}$, and we obtain
$$
a>
a \gcd(b_1-b_2, 1)^2=
\gcd(b_1,a)^2. 
\eqno\square
$$

\section{Level, Twist and $q$-Expansion  of a $j$-Map}
\label{sjmaps}
In this section we collect some properties of $j$-maps used in the sequel. %Notation introduced in this section will be widely used throughout the article. 

Given ${\gamma,\gamma'\in \GL_2^+(\Q)}$, we have ${j(\gamma z)=j(\gamma'z)}$ if and only if the matrices~$\gamma$ and~$\gamma'$ are \textsl{equivalent} in the sense of Section~\ref{srmat}. Combined with Proposition~\ref{ptri}, this gives the following.

%We will use the following description of $j$-maps.

\begin{proposition}
Let~$f$ be a non-constant $j$-map. Then there exist a unique positive number ${m\in\Q}$  
and a unique modulo~$1$ number ${\mu\in\Q}$  such that ${f(z)=j(mz+\mu)}$. \qed
\end{proposition}

Note that  ${m=\nlc(\gamma)}$ for any ${\gamma\in \GL_2^+(\Q)}$  such that ${f(z)=j(\gamma z)}$.

Setting ${q=e^{2\pi iz}}$ and ${\eps=e^{2\pi i\mu}}$, the map ${f(z)=j(mz+\mu)}$
admits the ``$q$-expansion'' 
\begin{equation}
\label{efqexp}
f(z)=\eps^{-1}q^{-m}+744+196884\eps q^{m}+21493760\eps^2q^{2m}+o(q^{2m}),
%864299970\eps^3q^{3m}+\ldots%20245856256\eps^4q^{4m}\ldots,
\end{equation}
where here and below we accept the following convention:
\begin{itemize}
%\item ``$\cdots$'' in a $q$-expansion means ``terms of $q$-degree higher than any of the preceding terms'';
\item
$O(q^\ell)$ means ``terms of $q$-degree~$\ell$ or higher'';
\item
$o(q^\ell)$ means ``terms of $q$-degree strictly higher than~$\ell$''.

\end{itemize} 

We call~$m$ and~$\eps$ the \textsl{level} and the \textsl{twist} of the non-constant $j$-map~$f$. For a constant $j$-map we set its level to be~$0$ and its twist undefined. 
The following property will be routinely used, usually without special reference:
\begin{equation}
\label{ef=glevtwist}
\text{two non-constant $j$-maps coincide if and only if their levels and twists coincide}. 
\end{equation}
We will denote in the sequel
%\begin{equation*}
%\label{eabc}
${A=196884}$ and ${B=21493760}$, %\quad C=864299970,%\quad D=20245856256,
%\end{equation*}
so that~\eqref{efqexp} reads 
\begin{equation}
\label{efqexpabc}
f(z)=\eps^{-1}q^{-m}+744+A\eps q^{m}+B\eps^2q^{2m}+O(q^{2m})
%C\eps^3q^{3m}+\ldots %D\eps^4q^{4m}\ldots.
\end{equation}

%Propositions~\ref{ptwo} and~\ref{pthree} translates into the language of $j$-maps as follows.

The following lemma will play an important role in Section~\ref{sini}. 

\begin{lemma}
\label{lthreejmaps}
Let $f_1$, $f_2$ and $f_3$ be pairwise distinct $j$-maps, not all constant. Then there exists ${\gamma\in \GL_2^+(\Q)}$ such that one of the maps ${f_1\circ\gamma}$, ${f_2\circ\gamma}$, ${f_3\circ\gamma}$ has level strictly bigger than the two others.
\end{lemma}

\paragraph{Proof}
If only one of the maps $f_k$ is non-constant then there is nothing to prove. If exactly two of them, say,~$f_1$ and~$f_2$, are non-constant, then Proposition~\ref{ptwo} implies the existence of ${\gamma\in \GL_2^+(\Q)}$ such that 
${f_1\circ\gamma}$ and ${f_2\circ\gamma}$ have distinct levels, and we are done. Finally, if all the three are non-constant, the result follows from Proposition~\ref{pthree}. \qed

\bigskip

We conclude this section with a linear independence property of non-constant $j$-maps.

\begin{lemma}
\label{lnlin}
Let $f,g$ be non-constant $j$-maps satisfying a non-trivial linear relation ${af+bg+c=0}$, where ${(a,b,c)\in \C^3}$ and ${(a,b,c)\ne (0,0,0)}$. Then ${f=g}$ and  ${a+b=c=0}$. 
\end{lemma}

\paragraph{Proof}
Any two non-constant $j$-maps parametrize the modular curve $Y_0(N)$ of a certain level~$N$; in other words, we have ${\Phi_N(f,g)=0}$, where ${\Phi_N(x,y)}$ is the $N$th modular polynomial. If we also have ${af+bg+c=0}$, then the polynomial $\Phi_N(x,y)$, being irreducible, must divide the linear polynomial ${ax+by+c}$. It follows that ${N=1}$, since ${\Phi_1(x,y)=x-y}$ is the only modular polynomial of degree~$1$. The result follows.\qed

\section{Initializing the Proof of the Main Lemma}
\label{sini}
In this section we start the proof of the Main Lemma. Thus, from now on, let  $f_1,f_2,f_3,g_1,g_2,g_3$ be $j$-maps, not all constant and satisfying
{%\small
\begin{equation}
\label{ecolli} 
\begin{vmatrix}
1&1&1\\
f_1&f_2&f_3\\
g_1&g_2&g_3
\end{vmatrix}=0.
\end{equation}}%
This can be rewritten as 
\begin{equation}
\label{edouble}
(f_1-f_2)(g_2-g_3)=(f_2-f_3)(g_1-g_2).
\end{equation}
If, say, ${f_1=f_2}$ then we find from~\eqref{edouble} that either ${f_2=f_3}$ in which case ${f_1=f_2=f_3}$, or ${g_1=g_2}$ in which case ${f_1=f_2}$ and ${g_1=g_2}$. Hence we may assume in the sequel that 
\begin{equation}
\label{edistinct}
\text{${f_1,f_2,f_3}$ are pairwise distinct, and so are ${g_1,g_2,g_3}$}.
\end{equation} 
We will show that under this assumption %either 
\begin{equation}
\label{erows}
f_k=g_k \qquad (k=1,2,3).
\end{equation}

Let $m_k,n_k$ be the levels of $f_k,g_k$, respectively, for ${k=1,2,3}$. If~$f_k$ and/or~$g_k$ is not constant, we denote the corresponding twists by ${\eps_k=e^{2\pi i\mu_k}}$ and/or ${\eta_k=e^{2\pi i\nu_k}}$, respectively.

\subsection{Some relations for the levels}
Since not all of our six maps are constant, we may assume that the three maps $f_k$ are  not all constant. Lemma~\ref{lthreejmaps} implies now that, after a suitable variable change, one of the numbers ${m_1,m_2,m_3}$ is strictly bigger than the others. After renumbering, we may assume that 
$$
m_1>m_2,m_3.
$$
We claim that 
\begin{equation}
\label{eng}
n_1>n_2,n_3
\end{equation} 
as well, and, moreover,
\begin{equation}
\label{edifference}
m_1-\max\{m_2,m_3\}= n_1-\max\{n_2,n_3\}. 
\end{equation}
Indeed, assume that, say,  ${n_2\ge n_1,n_3}$. Then the leading terms of the $q$-expansion on the left and on the right of~\eqref{edouble} are of the form ${cq^{-(m_1+n_2)}}$ and ${c'q^{-(\max\{m_2,m_3\}+n_2)}}$ with some non-zero~$c$ and~$c'$. (Precisely:
{%\small
$$
c=
\begin{cases}
\eps_1^{-1}\eta_2^{-1}, & n_2>n_3,\\
\eps_1^{-1}(\eta_2^{-1}-\eta_3^{-1}), & n_2=n_3>0,\\
\eps_1^{-1}(g_2-g_3), & n_2=n_3=0,
\end{cases}
$$}%
and it follows from~\eqref{edistinct} that ${c\ne 0}$; in a similar way one shows that ${c'\ne 0}$.) And this is impossible, because  ${m_1+n_2>\max\{m_2,m_3\}+n_2}$. This proves that ${n_1>n_2,n_3}$. In particular the three maps $g_k$ are also not all constant. Again comparing the leading terms of the $q$-expansion on the left and on the right of~\eqref{edouble}, we obtain~\eqref{edifference}. 

Swapping, if necessary, the functions $f_k$ and $g_k$, we may assume that 
\begin{equation}
\label{emonegenone}
m_1\ge n_1, 
\end{equation}
and after renumbering, we may assume that  
\begin{equation}
\label{emmm}
m_1> m_2\ge m_3.
\end{equation}
Equality~\eqref{edifference} now becomes 
\begin{equation}
\label{ediffbis}
m_1-m_2= n_1-\max\{n_2,n_3\}. 
\end{equation}

\subsection{One more lemma}

Here is a less obvious property, which will be used in the proof several times. 

\begin{lemma}
\label{lcross}
In the above set-up, we cannot have simultaneously ${f_2=g_3}$ and ${g_2=f_3}$. 
\end{lemma}
%This lemma is proved in Subsection~\ref{ssprlcross}. 

\paragraph{Proof}

If ${f_2=g_3}$ and ${g_2=f_3}$ then
{%\small
$$
0=
\begin{vmatrix}
1&1&1\\
f_1&f_2&f_3\\
g_1&f_3&f_2
\end{vmatrix}= (f_3-f_2)(f_1+g_1-f_2-f_3). 
$$}%
Since ${f_2\ne f_3}$, this implies 
\begin{equation}
\label{efgff}
f_1+g_1=f_2+f_3.
\end{equation}
We will see that this leads to a contradiction.

Observe first of all that ${m_2>0}$. Indeed, if ${m_2=0}$ then ${m_3=0}$ as well by~\eqref{emmm}. Hence both~$f_2$ and~$f_3$ are constant, and~\eqref{efgff} contradicts Lemma~\ref{lnlin}. 

Next, we have ${m_3>0}$ as well. Indeed, if~$f_3$ is constant, then, comparing the constant terms in~\eqref{efgff}, we find ${f_3=744}$, contradicting Lemma~\ref{ljne744etc}. 

Thus, we have 
${m_1\ge n_1> n_3=m_2\ge m_3>0}$. 
Comparing the $q$-expansions
\begin{align*}
f_1+g_1&=
\begin{cases}
\eps_1^{-1}q^{-m_1}+\eta_1^{-1}q^{-n_1}+O(1), & m_1>n_1,\\
(\eps_1^{-1}+\eta_1^{-1})q^{-m_1}+O(1), & m_1=n_1, \ \eps_1\ne-\eta_1,\\
1488+2B\eps_1^2q^{2m_1}+o(q^{2m_1}), & m_1=n_1, \  \eps_1=-\eta_1,
\end{cases}\\
f_2+f_3&=
\begin{cases}
\eps_2^{-1}q^{-m_2}+\eps_3^{-1}q^{-m_3}+O(1), & m_2>m_3, \\
(\eps_2^{-1}+\eps_3^{-1})q^{-m_2}+O(1), & m_2=m_3, \ \eps_2\ne-\eps_3,\\
1488+2B\eps_2^2q^{2m_2}+o(q^{2m_2}), & m_2=m_3, \ \eps_2=-\eps_3,\\
\end{cases}
\end{align*}
we immediately derive a contradiction.  \qed

\subsection{The determinant $\DD(q)$}

We will study in the sequel a slightly modified version of   the determinant from~\eqref{ecolli}:
{%\small
$$
\DD(q)=\begin{vmatrix}
1&1&1\\
q^{m_1}f_1&q^{m_1}f_2&q^{m_1}f_3\\
q^{n_1}g_1&q^{n_1}g_2&q^{n_1}g_3
\end{vmatrix}.
$$}%
The advantage is that it has no negative powers of~$q$. Equality~\eqref{ecolli} simply means that $\DD(q)$ vanishes as a formal power series in~$q$. It will be useful to write
{%\small
\begin{equation}
\label{e-744}
\DD(q)=\begin{vmatrix}
1&1&1\\
q^{m_1}(f_1-744)&q^{m_1}(f_2-744)&q^{m_1}(f_3-744)\\
q^{n_1}(g_1-744)&q^{n_1}(g_2-744)&q^{n_1}(g_3-744)
\end{vmatrix}
\end{equation}}%
This would allow us to eliminate the constant terms in the $q$-expansions of~$f_k$ and~$g_k$.

It will be convenient to use the notation 
\begin{equation}
\label{etilfg}
\tilf_k=
\begin{cases}
\eps_k^{-1},&m_k>0,\\
f_k-744,&m_k=0,
\end{cases},\quad
\tilg_k=
\begin{cases}
\eta_k^{-1},&n_k>0,\\
g_k-744,&n_k=0,
\end{cases}
\end{equation}
so that 
$$
q^{m_1}(f_k-744)=\tilf_kq^{m_1-m_k}+o(q^{m_1}), \quad q^{n_1}(g_k-744)=\tilg_kq^{n_1-n_k}+o(q^{n_1}).
$$
Lemma~\ref{ljne744etc} implies that 
\begin{equation}
\label{etilfgnezero}
\tilf_k, \tilg_k\ne 0 \qquad (k=1,2,3),
\end{equation}
which will be frequently used, usually without special references. 

\subsection{The four cases}

According to~\eqref{eng} and~\eqref{emmm}, there are four possible cases: 
{%\small
\begin{align*}
&m_2=m_3;\\
&m_2>m_3,\quad
n_2>n_3;\\
&m_2>m_3,\quad
n_2=n_3;\\
&m_2>m_3,\quad
n_3>n_2.
\end{align*}}%
They are treated in the four subsequent sections, respectively. %We consider them separately in the four subsequent sections. 
We will show that in the first two cases we have~\eqref{erows},  and the last two cases are impossible. 
The proofs in the four cases are similar in strategy but differ in technical details. 

Most of our arguments are nothing more than careful manipulations with $q$-expansions. Still, they are quite technical, and, to facilitate reading, we split proofs of each of the cases it into short logically complete steps.

\section{The Case ${m_2=m_3}$}
\label{smg=}

In this section we assume that 
$$
m_1>m_2=m_3.
$$ 
We want to prove that in this case we have
${f_k=g_k}$ for  ${k=1,2,3}$.  

Let us briefly describe the strategy of the proof. We already have~\eqref{eng}, and after renumbering we may assume that 
$$
n_1> n_2\ge n_3.
$$
Equality~\eqref{ediffbis} becomes now
\begin{equation}
\label{emmmmnn}
m_1-m_2=m_1-m_3=n_1-n_2.
\end{equation}
We start by proving that ${n_2=n_3}$, see Subsection~\ref{ssnsequal}. This been done, setting ${m_2=m_3=m}$ and ${n_2=n_3=n}$, we rewrite~\eqref{emmmmnn} as 
\begin{equation}
\label{emmnn}
m_1-m=n_1-n. 
\end{equation}
The next step is proving (see Subsection~\ref{ssmnequal}) that ${m_1=n_1}$.  In view of~\eqref{emmnn} this would imply that ${m=n}$ as well. In particular,~$f_k$ and~$g_k$ are of the same level for every ${k=1,2,3}$.  After this, we will be ready to prove that ${f_k=g_k}$ for ${k=1,2,3}$, see   Subsection~\ref{ssfgequal}.

\subsection{Proof of  ${n_2=n_3}$}
\label{ssnsequal}
In this subsection we prove  that
${n_2=n_3}$. 
Set 
$$
m_1-m_2=m_1-m_3=n_1-n_2=\lambda, \quad n_1-n_3=\lambda'\ge \lambda.
$$
We want to show that ${\lambda'=\lambda}$. 

Assume that ${\lambda'>\lambda}$. Then by (\ref{emonegenone}) all the~$m_k$ and~$n_k$ except perhaps~$n_3$ are positive. We consider separately the cases ${n_3=0}$ and ${n_3>0}$.

\subsubsection{The subcase ${n_3=0}$}
If ${n_3=0}$ then, using notation~\eqref{etilfg}, we write  ${\tilg_3=g_3-744}$ and 
{%\small
\begin{align*}
\DD(q)&=
\begin{vmatrix*}[l]
1&1&1\\
\eps_1^{-1}&\eps_2^{-1}q^\lambda&\eps_3^{-1}q^\lambda\\
\eta_1^{-1}&\eta_2^{-1}q^{\lambda}&\tilg_3q^{\lambda'}
\end{vmatrix*}+o(q^{n_1})\\
&=(\eps_1^{-1}\eta_2^{-1}-\eps_2^{-1}\eta_1^{-1}+\eps_3^{-1}\eta_1^{-1})q^{\lambda}+\eps_3^{-1}\eta_2^{-1}q^{2\lambda}+\eps_1^{-1}\tilg_3q^{\lambda'}+o(q^{n_1})+O(q^{\lambda+\lambda'}).
\end{align*}}%
The term with~$q^{\lambda'}$ can be eliminated only if ${\lambda'=2\lambda}$ and ${\eps_1^{-1}\tilg_3=\eps_3^{-1}\eta_2^{-1}}$, that is, 
${g_3=744+\eps_1\eps_3^{-1}\eta_2^{-1}}$,
contradicting Lemma~\ref{ljneroot}.% implies now that ${\eps_1\eps_3^{-1}\eta_2^{-1}=\pm1}$. Hence ${g_3=744\pm1}$, contradicting Lemma~\ref{ljne744etc}. 

\subsubsection{The subcase ${n_3>0}$}
If ${n_3>0}$ then 
{%\small
\begin{align*}
\DD(q)&=
\begin{vmatrix*}[l]
1&1&1\\
\eps_1^{-1}&\eps_2^{-1}q^\lambda&\eps_3^{-1}q^\lambda\\
\eta_1^{-1}&\eta_2^{-1}q^{\lambda}&\eta_3^{-1}q^{\lambda'}+A\eta_3q^{n_1+n_3}
\end{vmatrix*}+o(q^{n_1+n_3})\\
&=(\eps_1^{-1}\eta_2^{-1}-\eps_2^{-1}\eta_1^{-1}+\eps_3^{-1}\eta_1^{-1})q^{\lambda}-\eps_3^{-1}\eta_2^{-1}q^{2\lambda}-\eps_1^{-1}\eta_3^{-1}q^{\lambda'}+\eps_2^{-1}\eta_3^{-1}q^{\lambda+\lambda'}\\
&\hphantom{=}-A\eps_1^{-1}\eta_3q^{n_1+n_3}+o(q^{n_1+n_3}).
\end{align*}}%
As $n_1+n_3> \lambda'$, the term with $q^{n_1+n_3}$ can be eliminated only if either 
$$
\lambda<\lambda'<2\lambda=n_1+n_3<\lambda+\lambda', \quad \eps_3^{-1}\eta_2^{-1}=-A\eps_1^{-1}\eta_3,
$$
which is impossible because~$A$ is not a root of unity, or 
$$
\lambda<\lambda',2\lambda<n_1+n_3=\lambda+\lambda', \quad \eps_2^{-1}\eta_3^{-1}=A\eps_1^{-1}\eta_3,
$$
which is again impossible by the same reason.

\subsubsection{Conclusion}

Thus, we have proved that ${n_2=n_3}$. Setting ${m=m_2=m_3}$ and ${n=n_2=n_3}$, we can summarize our knowledge as follows:
\begin{equation*}
%\label{emnmnmn}
m_1>m_2=m_3=m; \quad n_1>n_2=n_3=n; \quad m_1-m=n_1-n=\lambda>0; \quad  m_1-n_1=m-n\ge 0. 
\end{equation*}
Together with~\eqref{edistinct} this implies that 
\begin{equation}
\label{eeneznz}
%\text{$\eps_2\ne \eps_3$ if ${m>0}$ and $\eta_2\ne\eta_3$ if ${n>0}$}. 
\tilf_2\ne\tilf_3, \quad \tilg_2\ne\tilg_3. 
\end{equation}

\subsection{Proof of ${m_1=n_1}$}
\label{ssmnequal}
Now we want to prove that 
\begin{equation}
\label{em=n}
m_1=n_1.
\end{equation} 
%By~\eqref{emnmnmn}, this is equivalent to  ${m=n}$. 
Thus, assume that ${m_1>n_1}$, in which case we also have ${m>n}$. We consider separately the subcases ${n>0}$ and ${n=0}$.

\subsubsection{The subcase ${n>0}$.}
%\label{sssngz}
If ${n>0}$ then %, using~\eqref{efqexpabc}, % and denoting ${\ell=\min\{m_1+m,2n_1, n_1+2n\}}$, 
%we write~\eqref{e-744} as
{%\small
\begin{align*}
\DD(q)&=
\begin{vmatrix*}[l]
1&1&1\\
\eps_1^{-1}&\eps_2^{-1}q^\lambda&\eps_3^{-1}q^\lambda\\
\eta_1^{-1}&\eta_2^{-1}q^{\lambda}+A\eta_2q^{n_1+n}&\eta_3^{-1}q^{\lambda}+A\eta_3q^{n_1+n}
\end{vmatrix*}+o(q^{n_1+n})\\ 
&=
\begin{vmatrix}
\eps_1^{-1}&\eps_2^{-1}-\eps_3^{-1}\\
\eta_1^{-1}&\eta_2^{-1}-\eta_3^{-1}
\end{vmatrix}q^{\lambda}+ 
\begin{vmatrix}
\eps_2^{-1}&\eps_3^{-1}\\
\eta_2^{-1}&\eta_3^{-1}
\end{vmatrix}q^{2\lambda}+ A\eps_1^{-1}(\eta_2-\eta_3)q^{n_1+n}+o(q^{n_1+n})+o(q^{2\lambda}).
\end{align*}}%
%Now if ${2\lambda\ne n_1+n}$ then the term involving $q^{n_1+n}$ is non-zero by~\eqref{eeneznz}, % which implies ${\DD(q)\ne 0}$, a contradiction. And if 
Here the coefficient of $q^\lambda$ must vanish. If ${2\lambda > n_1+n}$ then that of $q^{n_1+n}$ must vanish too; but that would contradict~\eqref{eeneznz}. If ${2\lambda < n_1+n}$ then the coefficient of $q^{2\lambda}$ must vanish and then that of $q^{n_1+n}$. It follows that ${2\lambda= n_1+n}$ and 
\begin{equation}
\label{egal}
\begin{vmatrix}
\eps_2^{-1}&\eps_3^{-1}\\
\eta_2^{-1}&\eta_3^{-1}
\end{vmatrix}=A\eps_1^{-1}(\eta_3-\eta_2).
\end{equation}
As noted, both sides of~\eqref{egal} are non-zero. Since the left-hand side is a sum of~$2$ roots of unity, Lemma~\ref{lsumroots} implies that ${196884=|A|\le 2}$, a contradiction. 
This completes the proof of~\eqref{em=n} in the case ${n>0}$.

\subsubsection{The subcase ${n=0}$.}
If ${n=0}$ then~$g_2$ and~$g_3$ are distinct constants, and the other functions are non-constant. Also, we have ${\lambda=n_1}$, and so
\begin{equation}
\label{emonenone}
m_1=m+n_1. 
\end{equation} 
Now, using notation~\eqref{etilfg}, we obtain %denote ${\ell=\min\{2m_1,m_1+2m,3n_1\}}$ and write~\eqref{e-744} as
{%\small
\begin{align*}
\DD(q)&=
\begin{vmatrix*}[l]
1&1&1\\
\eps_1^{-1}&\eps_2^{-1}q^{n_1}+A\eps_2q^{m_1+m}&\eps_3^{-1}q^{n_1}+A\eps_3q^{m_1+m}\\
\eta_1^{-1}+A\eta_1q^{2n_1}&\tilg_2q^{n_1}&\tilg_3q^{n_1}
\end{vmatrix*}+o(q^{m_1+m})+o(q^{2n_1})\\ 
&=
\begin{vmatrix}
\eps_1^{-1}&\eps_2^{-1}-\eps_3^{-1}\\
\eta_1^{-1}&\tilg_2-\tilg_3
\end{vmatrix}q^{n_1}+ 
\begin{vmatrix}
\eps_2^{-1}&\eps_3^{-1}\\
\tilg_2&\tilg_3
\end{vmatrix}q^{2n_1}+ A\eta_1^{-1}(\eps_3-\eps_2)q^{m_1+m}+o(q^{m_1+m})+o(q^{2n_1}).
\end{align*}}%
%Now if ${2n_1\ne m_1+m}$ then the term involving $q^{m_1+m}$ is non-zero by~\eqref{eeneznz}, and  ${\DD(q)\ne 0}$, a contradiction. 
%We are left with the case 
As ${\eps_3 \neq \eps_2}$, the coefficient of $q^{m_1+m}$ is non-zero; by Lemma~\ref{lrootofone} so is the coefficient of $q^{2n_1}$.  This shows that ${2n_1= m_1+m}$. Together with~\eqref{emonenone} %and~\eqref{emknk} 
this implies ${m_1=3m}$ and ${n_1=2m}$; rescaling~$z$, we may assume
$$
m=1, \quad n_1=2, \quad m_1=3. 
$$ 
Hence
{%\small
\begin{align*}
\DD(q)&=
\begin{vmatrix*}[l]
1&1&1\\
\eps_1^{-1}&\eps_2^{-1}q^{2}+A\eps_2q^{4}+B\eps_2^2q^5&\eps_3^{-1}q^{2}+A\eps_3q^{4}+B\eps_3^2q^5\\
\eta_1^{-1}+A\eta_1q^{4}&\tilg_2q^{2}&\tilg_3q^{2}
\end{vmatrix*}+O(q^{6})\\ 
&=
\begin{vmatrix}
\eps_1^{-1}&\eps_2^{-1}-\eps_3^{-1}\\
\eta_1^{-1}&\tilg_2-\tilg_3
\end{vmatrix}q^{2}+ 
\left(\begin{vmatrix}
\eps_2^{-1}&\eps_3^{-1}\\
\tilg_2&\tilg_3
\end{vmatrix}+ A\eta_1^{-1}(\eps_3-\eps_2)\right)q^{4}+B\eta_1^{-1}(\eps_3^2-\eps_2^2)q^5+O(q^6).
\end{align*}}%
Equating to~$0$ the coefficient of~$q^5$, we obtain ${\eps_3=\pm\eps_2}$, and~\eqref{eeneznz} implies that ${\eps_3=-\eps_2}$. 
Using this, and equating to~$0$ the coefficients of~$q^2$ and~$q^4$, we obtain
$$
\eps_1^{-1}(\tilg_2-\tilg_3)=2\eps_2^{-1}\eta_1^{-1}, \quad \eps_2^{-1}(\tilg_2+\tilg_3)=2A\eta_1^{-1}\eps_2. 
$$
from which we deduce ${g_2=\tilg_2+744=\eps_1\eps_2^{-1}\eta_1^{-1}+A\eta_1^{-1}\eps_2^2+744}$.

Now Lemma~\ref{ljnetworoots} implies that ${g_2\in \Z}$, from which we deduce, using Lemma~\ref{lquadratic}, that both roots of unity ${\eps_1\eps_2^{-1}\eta_1^{-1}}$ and ${\eta_1^{-1}\eps_2^2}$ must be $\pm1$. Hence~$g_2$ is one of the four numbers ${744\pm1\pm A}$, contradicting Lemma~\ref{ljne744etc}.

\subsection{Proof of ${f_k=g_k}$ for ${k=1,2,3}$}
\label{ssfgequal}
In the previous subsection we proved that 
\begin{equation}
\label{emngmmnn}
m_1=n_1> m=n.
\end{equation}
We want to prove now that
\begin{equation}
\label{ef=g}
f_k=g_k \qquad (k=1,2,3). 
\end{equation} 
We again distinguish the subcases ${m=n>0}$ and  ${m=n=0}$. As before, we set ${\lambda=m_1-m=n_1-n}$. 

\subsubsection{The subcase ${m=n>0}$.}
If ${m=n>0}$ then %,  denoting  %${A=196884}$ and ${\ell=\min\{2\lambda+2m,\lambda+3m\}}$, we write~\eqref{e-744} as
{%\small
\begin{align}
\DD(q)&=
\begin{vmatrix*}[l]
1&1&1\\
\eps_1^{-1}&\eps_2^{-1}q^\lambda+A\eps_2q^{\lambda+2m}&\eps_3^{-1}q^\lambda+A\eps_3q^{\lambda+2m}\\
\eta_1^{-1}&\eta_2^{-1}q^{\lambda}+A\eta_2q^{\lambda+2m}&\eta_3^{-1}q^{\lambda}+A\eta_3q^{\lambda+2m}
\end{vmatrix*}+o(q^{\lambda+2m})\nonumber\\ 
\label{edq}
&=
\begin{vmatrix}
\eps_1^{-1}&\eps_2^{-1}-\eps_3^{-1}\\
\eta_1^{-1}&\eta_2^{-1}-\eta_3^{-1}
\end{vmatrix}q^{\lambda}+ 
\begin{vmatrix}
\eps_2^{-1}&\eps_3^{-1}\\
\eta_2^{-1}&\eta_3^{-1}
\end{vmatrix}q^{2\lambda}+ A\begin{vmatrix}
\eps_1^{-1}&\eps_2-\eps_3\\
\eta_1^{-1}&\eta_2-\eta_3
\end{vmatrix}q^{\lambda+2m}+o(q^{\lambda+2m}).
\end{align}}%
%Now we argue differently depending on whether ${\lambda\ne 2m}$ or ${\lambda=2m}$.
%\paragraph{\underline{The subcase ${\lambda\ne 2m}$}} 
%If ${\lambda\ne 2m}$ then 
This implies the equations
{%\small
\begin{equation}
\label{eequs}
\begin{vmatrix}
\eps_1^{-1}&\eps_2^{-1}-\eps_3^{-1}\\
\eta_1^{-1}&\eta_2^{-1}-\eta_3^{-1}
\end{vmatrix}=0, \quad
\begin{vmatrix}
\eps_2^{-1}&\eps_3^{-1}\\
\eta_2^{-1}&\eta_3^{-1}
\end{vmatrix}=0,\quad\begin{vmatrix}
\eps_1^{-1}&\eps_2-\eps_3\\
\eta_1^{-1}&\eta_2-\eta_3
\end{vmatrix}=0
\end{equation}}%
if ${2\lambda\ne \lambda+2m}$, and the equations
{%\small
\begin{gather}
\begin{vmatrix}
\eps_1^{-1}&\eps_2^{-1}-\eps_3^{-1}\\
\eta_1^{-1}&\eta_2^{-1}-\eta_3^{-1}
\end{vmatrix}=0, \nonumber\\
\label{eaab}
\begin{vmatrix}
\eps_2^{-1}&\eps_3^{-1}\\
\eta_2^{-1}&\eta_3^{-1}
\end{vmatrix}=-A\begin{vmatrix}
\eps_1^{-1}&\eps_2-\eps_3\\
\eta_1^{-1}&\eta_2-\eta_3
\end{vmatrix}
\end{gather}}%
if ${2\lambda= \lambda+2m}$. 
If both sides of~\eqref{eaab} are non-zero, then Lemma~\ref{lsumroots} implies ${196884=|A|\le 2}$,   a contradiction. Hence in any case we have~\eqref{eequs}.

Resolving the first two equations from~\eqref{eequs} in ${\eta_1^{-1},\eta_2^{-1},\eta_3^{-1}}$ and using~\eqref{eeneznz}, we obtain
$$
(\eta_1,\eta_2,\eta_3)=\theta(\eps_1,\eps_2,\eps_3)
$$
for some ${\theta\in \C}$. Substituting this to the third equation in~\eqref{eequs} and again using~\eqref{eeneznz}, we find ${\theta= \pm1}$. 
If ${\theta=-1}$ then
{%\small
\begin{align*}
\DD(q)&=
\begin{vmatrix*}[l]
\hphantom{-}1&\hphantom{-}1&\hphantom{-}1\\
\hphantom{-}\eps_1^{-1}+A\eps_1 q^{2\lambda+2m}&\hphantom{-}\eps_2^{-1}q^\lambda+A\eps_2q^{\lambda+2m}+B\eps_2^2q^{\lambda+3m}&\hphantom{-}\eps_3^{-1}q^\lambda+A\eps_3q^{\lambda+2m}+B\eps_3^2q^{\lambda+3m}\\
-\eps_1^{-1}-A\eps_1 q^{2\lambda+2m}&-\eps_2^{-1}q^\lambda-A\eps_2q^{\lambda+2m}+B\eps_2^2q^{\lambda+3m}&-\eps_3^{-1}q^\lambda-A\eps_3q^{\lambda+2m}+B\eps_3^2q^{\lambda+3m}
\end{vmatrix*}\\
&\hphantom{=}+ o(q^{\lambda+3m})\\ 
&=
\begin{vmatrix*}[l]
1&1&1\\
\eps_1^{-1}+A\eps_1 q^{2\lambda+2m}&\eps_2^{-1}q^\lambda+A\eps_2q^{\lambda+2m}+B\eps_2^2q^{\lambda+3m}&\eps_3^{-1}q^\lambda+A\eps_3q^{\lambda+2m}+B\eps_3^2q^{\lambda+3m}\\
&2B\eps_2^2q^{\lambda+3m}&2B\eps_3^2q^{\lambda+3m}
\end{vmatrix*}\\
&\hphantom{=}+o(q^{\lambda+3m})\\ 
&=2B\eps_1^{-1}(\eps_2^2-\eps_3^2)q^{\lambda+3m}+ o(q^{\lambda+3m}),
\end{align*}}%
which gives ${\eps_2=\pm\eps_3}$, and ${\eps_2=-\eps_3}$ by~\eqref{eeneznz}. 
Thus, we have ${\eps_2=\eta_3=-\eps_3=-\eta_2}$, which implies that ${f_2=g_3}$ and ${g_2=f_3}$, contradicting Lemma~\ref{lcross}. 

The only remaining option is ${\theta=1}$, which, together with~\eqref{emngmmnn},   proves~\eqref{ef=g}.

\subsubsection{The subcase ${m=n=0}$}
This case can be easily settled using Lemma~\ref{lnlin}. Indeed, in the case ${m=n=0}$ the functions $f_1,g_1$ are non-constant,  ${f_2,f_3,g_2,g_3}$ are constant, and 
{%\small
\begin{equation*}
0=\begin{vmatrix}
1&1&1\\
f_1&f_2&f_3\\
g_1&g_2&g_3
\end{vmatrix}=(g_2-g_3)f_1-(f_2-f_3)g_1+ 
\begin{vmatrix}
f_2&f_3\\
g_2&g_3
\end{vmatrix}
\end{equation*}}%
is  a non-trivial  linear relation for $f_1,g_1$ (recall that ${f_2\ne f_3}$ and ${g_2\ne g_3}$ by~\eqref{edistinct}.) By Lemma~\ref{lnlin}
{%\small
$$
f_1=g_1, \quad f_2-f_3=g_2-g_3, \quad \begin{vmatrix}
f_2&f_3\\
g_2&g_3
\end{vmatrix}=0.
$$}%
From the last two equations one easily deduces that ${f_2=g_2}$ and ${f_3=g_3}$, proving~\eqref{ef=g}.

\section{The Case ${m_2>m_3}$,\quad  ${n_2> n_3}$}
In this section we assume that
\begin{equation}
\label{emnggg}
m_1>m_2>m_3, \quad n_1>n_2> n_3.
\end{equation}
As in the previous section, we will prove that in this case ${f_k=g_k}$ for ${k=1,2,3}$.

The strategy of the proof is similar to that of the previous section. 
Equality~\eqref{ediffbis} now reads 
\begin{equation}
\label{emone-mtwo=none-ntwo}
m_1-m_2=n_1-n_2.
\end{equation}
We start with proving that 
\begin{equation}
\label{emone-mthree=none-nthree}
m_1-m_3=n_1-n_3,
\end{equation}
see Subsection~\ref{ssm-m=n-n}. Next to this, we prove, in Subsection~\ref{ssmone=none}, that ${m_1=n_1}$. Since, by this time, we will already know~\eqref{emone-mtwo=none-ntwo} and~\eqref{emone-mthree=none-nthree}, this will imply that ${m_k=n_k}$ for every ${k=1,2,3}$. After this, we prove that ${f_k=g_k}$ for  ${k=1,2,3}$ in Subsection~\ref{ssthetaone}. 

\bigskip

We set ${m_1-m_2=n_1-n_2=\lambda}$. We also have ${m_1\ge n_1}$ by~\eqref{emonegenone}. Let us collect our knowledge:
$$
m_1>m_2>m_3; \quad n_1>n_2> n_3; \quad m_1-m_2=n_1-n_2=\lambda>0; \quad m_1-n_1=m_2-n_2\ge0. 
$$

\subsection{Proof of ${m_1-m_3=n_1-n_3}$}
\label{ssm-m=n-n}

Now let us prove that ${m_1-m_3=n_1-n_3}$. Using notation~\eqref{etilfg}, we write
{%\small
\begin{align*}
\DD(q)&=
\begin{vmatrix*}[l]
1&1&1\\
\eps_1^{-1}&\eps_2^{-1}q^{\lambda}&\tilf_3q^{m_1-m_3}\\
\eta_1^{-1}&\eta_2^{-1}q^{\lambda}&\tilg_3q^{n_1-n_3}
\end{vmatrix*}+o(q^{n_1})\\
&= \begin{vmatrix}
\eps_1^{-1}&\eps_2^{-1}\\
\eta_1^{-1}&\eta_2^{-1}
\end{vmatrix}q^\lambda+\tilf_3\eta_1^{-1}q^{m_1-m_3}-\eps_1^{-1}\tilg_3q^{n_1-n_3}+o(q^{m_1-m_3})+o(q^{n_1-n_3}).
\end{align*}}%
If ${m_1-m_3\ne n_1-n_3}$ then we have one of the following options:
$$
\lambda<m_1-m_3<n_1-n_3; \quad \lambda <n_1-n_3<m_1-m_3. 
$$
In the first case $q^{m_1-m_3}$ cannot be eliminated, and in the second case  $q^{n_1-n_3}$ cannot be eliminated. This proves that ${m_1-m_3=n_1-n_3}$. 

We set ${m_1-m_3=n_1-n_3=\lambda'}$.  Thus, we have 
{%\small
\begin{align}
&m_1>m_2>m_3; \quad n_1>n_2> n_3; \nonumber\\
&m_1-m_2=n_1-n_2=\lambda>0; \quad m_1-m_3=n_1-n_3=\lambda'>\lambda>0; \nonumber\\
\label{edelt}
&m_1-n_1=m_2-n_2=m_3-n_3\ge0. 
\end{align}}%
In addition to this, from 
{%\small
\begin{equation*}
\DD(q)=
\begin{vmatrix*}[l]
1&1&1\\
\eps_1^{-1}&\eps_2^{-1}q^{\lambda}&\tilf_3q^{\lambda'}\\
\eta_1^{-1}&\eta_2^{-1}q^{\lambda}&\tilg_3q^{\lambda'}
\end{vmatrix*}+o(q^{n_1})
= \begin{vmatrix}
\eps_1^{-1}&\eps_2^{-1}\\
\eta_1^{-1}&\eta_2^{-1}
\end{vmatrix}q^\lambda-
\begin{vmatrix}
\eps_1^{-1}&\tilf_3\\
\eta_1^{-1}&\tilg_3
\end{vmatrix}q^{\lambda'}
+o(q^{\lambda'}),
\end{equation*}}%
we deduce that 
{%\small
\begin{equation}
\label{edetszero}
\begin{vmatrix}
\eps_1^{-1}&\eps_2^{-1}\\
\eta_1^{-1}&\eta_2^{-1}
\end{vmatrix}=
\begin{vmatrix}
\eps_1^{-1}&\tilf_3\\
\eta_1^{-1}&\tilg_3
\end{vmatrix}=0,
\end{equation}}%
which means that 
\begin{equation}
\label{epropor}
(\eta_1^{-1},\eta_2^{-1},\tilg_3)=\theta (\eps_1^{-1},\eps_2^{-1},\tilf_3)
\end{equation}
with some root of unity~$\theta$. 

\subsection{Proof of  ${m_1=n_1}$}
\label{ssmone=none}
In this subsection we show that ${m_1=n_1}$. Thus, assume that 
\begin{equation}
\label{emonegnone}
m_1>n_1, 
\end{equation}
in which case we  also have 
\begin{equation}
\label{emsgns}
m_2>n_2, \quad m_3>n_3.
\end{equation}
We should also have 
\begin{equation}
\label{enthreepos}
n_3>0.
\end{equation} 
Indeed, if ${m_3>n_3=0}$ then the second equation in~\eqref{edetszero} reads 
${g_3=744+\eps_1\eps_3^{-1}\eta_1^{-1}}$, which is impossible by Lemma~\ref{ljneroot}.  %. Lemma~\ref{ljnecyc} implies that ${\eps_1\eps_3^{-1}\eta_1^{-1}=\pm1 }$ and ${g_3\in\{743,745\}}$, contradicting Lemma~\ref{ljne744}. 

Using~\eqref{epropor},~\eqref{emonegnone},~\eqref{emsgns},~\eqref{enthreepos}, we obtain
{%\small
\begin{align*}
\DD(q)&=
\begin{vmatrix*}[l]
1&1&1\\
\eps_1^{-1}&\eps_2^{-1}q^{\lambda}&\eps_3^{-1}q^{\lambda'}\\
\eta_1^{-1}&\eta_2^{-1}q^{\lambda}&\eta_3^{-1}q^{\lambda'}+A\eta_3q^{n_1+n_3}
\end{vmatrix*}+o(q^{n_1+n_3})\\
&= \begin{vmatrix*}[l]
1&1&1\\
\eps_1^{-1}&\eps_2^{-1}q^{\lambda}&\eps_3^{-1}q^{\lambda'}\\
0&0&A\eta_3q^{n_1+n_3}
\end{vmatrix*}+o(q^{n_1+n_3})\\
&=-A\eps_1^{-1}\eta_3q^{n_1+n_3}+o(q^{n_1+n_3}),
\end{align*}}%
a contradiction.

This proves that 
\begin{equation}
\label{emk=nkallk}
m_k=n_k \qquad (k=1,2,3).
\end{equation}

\subsection{Proof of ${f_k=g_k}$ for ${k=1,2,3}$}
\label{ssthetaone}
To prove that ${f_k=g_k}$ for ${k=1,2,3}$, we only need  to show that
$$
\theta=1,
$$
where~$\theta$ is from~\eqref{epropor}. If ${m_3=n_3=0}$ then, rewriting the equality ${\tilg_3=\theta\tilf_3}$ as ${(g_3-744)=\theta(f_3-744)}$, we deduce ${\theta=1}$  from  Lemma~\ref{lrootofone}. 

Now assume that  ${m_3=n_3>0}$. In this case
%\subsubsection{The subcase ${m_3=n_3>0}$}
%We prove first that ${\theta=\pm1}$. 
%If ${m_3=n_3>0}$ then
{%\small
\begin{align*}
\DD(q)&=
\begin{vmatrix*}[l]
1&1&1\\
\eps_1^{-1}&\eps_2^{-1}q^{\lambda}&\eps_3^{-1}q^{\lambda'}+A\eps_3q^{m_1+m_3}\\
\eta_1^{-1}&\eta_2^{-1}q^{\lambda}&\eta_3^{-1}q^{\lambda'}+A\eta_3q^{m_1+m_3}
\end{vmatrix*}+o(q^{m_1+m_3})\\
&=\begin{vmatrix*}[l]
1&1&1\\
\eps_1^{-1}&\eps_2^{-1}q^{\lambda}&\eps_3^{-1}q^{\lambda'}+A\eps_3q^{m_1+m_3}\\
0&0&A\eps_3(\theta^{-1}-\theta)q^{m_1+m_3}
\end{vmatrix*}+o(q^{m_1+m_3})\\
&=
-A\eps_1^{-1}\eps_3(\theta^{-1}-\theta)q^{m_1+m_3}
+o(q^{m_1+m_3}),
\end{align*}}%
which implies ${\theta=\pm1}$. If ${\theta=-1}$ then
{%\small
\begin{align*}
\DD(q)&=
\begin{vmatrix*}[l]
\hphantom{-}1&\hphantom{-}1&\hphantom{-}1\\
\hphantom{-}\eps_1^{-1}+A\eps_1q^{2m_1}&\hphantom{-}\eps_2^{-1}q^{\lambda}+A\eps_2q^{m_1+m_2}&\hphantom{-}\eps_3^{-1}q^{\lambda'}+A\eps_3q^{m_1+m_3}+B\eps_3^2q^{m_1+2m_3}\\
-\eps_1^{-1}-A\eps_1q^{2m_1}&-\eps_2^{-1}q^{\lambda}-A\eps_2q^{m_1+m_2}&-\eps_3^{-1}q^{\lambda'}-A\eps_3q^{m_1+m_3}+B\eps_3^2q^{m_1+2m_3}
\end{vmatrix*}+o(q^{m_1+2m_3})\\
&=\begin{vmatrix*}[l]
1&1&1\\
\eps_1^{-1}+A\eps_1q^{2m_1}&\eps_2^{-1}q^{\lambda}+A\eps_2q^{m_1+m_2}&\eps_3^{-1}q^{\lambda'}+A\eps_3q^{m_1+m_3}\\
0&0&2B\eps_3^2q^{m_1+2m_3}
\end{vmatrix*}+o(q^{m_1+2m_3})\\
&=
-2B\eps_1^{-1}\eps_3^2q^{m_1+2m_3}
+o(q^{m_1+2m_3}),
\end{align*}}%
a contradiction. %This proves that ${\theta=1}$ when ${m_3=n_3>0}$. 

Thus, in any case we have ${\theta=1}$ in~\eqref{epropor}. Together with~\eqref{emk=nkallk} this proves that ${f_k=g_k}$ for ${k=1,2,3}$.

\section{The Case ${m_2>m_3}$,\quad ${n_2=n_3}$}

In this section we assume that 
\begin{equation}
\label{emggng=}
m_1>m_2>m_3,\quad n_1>n_2=n_3,
\end{equation}
and will show that this is impossible. %We treat separately two subcases: ${n_2=n_3>0}$ and ${n_2=n_3=0}$. 

Relation~\eqref{ediffbis} now becomes ${m_1-m_2=n_1-n_2=n_1-n_3}$. We set
\begin{equation}
\label{ethreelams}
m_1-m_2=n_1-n_2=n_1-n_3=\lambda.
\end{equation}

Fist of all, let us rule out the case ${n_2=n_3=0}$. In this case ${n_1=\lambda<m_1-m_3}$. Using notation~\eqref{etilfg}, we write in this case
{%\small
\begin{align*}
\DD(q)&=
\begin{vmatrix*}[l]
1&1&1\\
\eps_1^{-1}&\eps_2^{-1}q^\lambda&0\\
\eta_1^{-1}&\tilg_2q^\lambda&\tilg_3q^\lambda
\end{vmatrix*}+o(q^\lambda)=(\eps_1^{-1}\tilg_2-\eps_1^{-1}\tilg_3-\eps_2^{-1}\eta_1^{-1})q^\lambda+o(q^\lambda). 
\end{align*}}%
We obtain 
${\eps_1^{-1}\tilg_2-\eps_1^{-1}\tilg_3-\eps_2^{-1}\eta_1^{-1}=0}$, which contradicts Lemma~\ref{lcute}. 

Thus, we may assume in the sequel that 
\begin{equation}
\label{eng=g}
n_2= n_3>0. 
\end{equation}
Since ${n_2=n_3}$, we have
\begin{equation}
\label{ezetasdif}
\eta_2\ne\eta_3,
\end{equation}
which will be systematically used, sometimes without special reference. 

Our principal objective will be to show that ${m_3=m_1-2\lambda}$ and ${n_1=m_1-\lambda/2}$. 
The first of these two relations is proved already in Subsection~\ref{sslmmnn}. The second one is more delicate and will be established in Subsection~\ref{ssmmnn<l}, after some preparatory work done in the previous subsections. On the way, we will also prove certain inequalities relating the numbers~$m_k$,~$n_k$ and~$\lambda$, and certain relations for the twists. After all this is done, obtaining a contradiction will be relatively easy, see Subsection~\ref{ssconcl}.

\subsection{Proof of ${2\lambda=m_1-m_3\le n_1+n_2}$}
\label{sslmmnn}
Using notation~\eqref{etilfg}, we write
{%\small
\begin{align}
\DD(q)&=
\begin{vmatrix*}[l]
1&1&1\\
\eps_1^{-1}&\eps_2^{-1}q^{\lambda}&\tilf_3q^{m_1-m_3}\\
\eta_1^{-1}&\eta_2^{-1}q^{\lambda}+A\eta_2q^{n_1+n_2}&\eta_3^{-1}q^{\lambda}+A\eta_3q^{n_1+n_2}
\end{vmatrix*}+o(q^{m_1})+o(q^{n_1+n_2})\nonumber\\
&=(\eps_1^{-1}\eta_2^{-1}-\eps_1^{-1}\eta_3^{-1}-\eps_2^{-1}\eta_1^{-1})q^\lambda+\eps_2^{-1}\eta_3^{-1}q^{2\lambda}+\eta_1^{-1}\tilf_3q^{m_1-m_3}+A\eps_1^{-1}(\eta_2-\eta_3)q^{n_1+n_2}\nonumber\\
\label{eterms}
&\hphantom{=}+o(q^{m_1-m_3})+o(q^{n_1+n_2}).
\end{align}}%
First of all, this gives 
\begin{equation}
\label{esrut}
\eps_1^{-1}\eta_2^{-1}-\eps_1^{-1}\eta_3^{-1}-\eps_2^{-1}\eta_1^{-1}=0.
\end{equation}
A sum of~$3$ roots of unity can vanish only if they are proportional to the~$3$ distinct cubic roots of unity. In particular, 
\begin{equation}
\label{esixth}
\text{$\eta_2/\eta_3$ is a primitive $6$th root of unity}. 
\end{equation}

\paragraph{\underline{We have ${m_1-m_3\ge 2\lambda}$}}

Indeed, if ${2\lambda>m_1-m_3}$ then we must have 
\begin{equation}\label{newlabel}
m_1-m_3=n_1+n_2, \quad \eta_1^{-1}\tilf_3=-A\eps_1^{-1}(\eta_2-\eta_3).
\end{equation} 
If ${m_3>0}$ this gives ${\eta_1^{-1}\eps_3^{-1}=-A\eps_1^{-1}(\eta_2-\eta_3)}$ which is impossible because~$A$ does not divide a root of unity. And if ${m_3=0}$ then ${f_3=744-A\eps_1^{-1}\eta_1(\eta_2-\eta_3)}$. 
Lemma~\ref{ljnetworoots} now implies that ${f_3\in \Z}$, and we obtain ${f_3\in \{744\pm196884,744\pm2\cdot196884\}}$ contradicting Lemma~\ref{ljne744etc}.

\paragraph{\underline{We have ${m_1-m_3\le 2\lambda}$}}

Indeed, if ${2\lambda<m_1-m_3}$ then the term with $q^{2\lambda}$ cancels either against a term in $o(q^{n_1+n_2})$ or against the term with $q^{n_1+n_2}$. In the first situation the terms with $q^{m_1-m_3}$ and $q^{n_1+n_2}$ must cancel each other, and we are back to (\ref{newlabel}). In the second situation we must have 
$$
2\lambda=n_1+n_2, \quad \eps_2^{-1}\eta_3^{-1}=-A\eps_1^{-1}(\eta_2-\eta_3),
$$ 
which is impossible because ${A=196884}$ does not divide a root of unity.

\bigskip

Thus, we proved that ${m_1-m_3=2\lambda}$. 

\paragraph{\underline{We have ${n_1+n_2\ge 2\lambda}$}}
Indeed, if ${n_1+n_2< 2\lambda=m_1-m_3}$ then the non-zero term ${A\eps_1^{-1}(\eta_2-\eta_3)q^{n_1+n_2}}$ cannot be eliminated. (It is non-zero because of~\eqref{ezetasdif}.) 

\bigskip

Thus, we proved that 
\begin{equation}
\label{einequ}
2\lambda=m_1-m_3\le n_1+n_2.
\end{equation}

\subsection{Proof of ${n_1+n_2>2\lambda}$}

We want to show now that the inequality in~\eqref{einequ} is strict. Thus, assume the contrary, that is,  
\begin{equation}
\label{ealleq}
2\lambda=m_1-m_3= n_1+n_2.
\end{equation}
Then~\eqref{eterms} implies that 
\begin{equation}
\label{esumofthree}
\eps_2^{-1}\eta_3^{-1}+\eta_1^{-1}\tilf_3+A\eps_1^{-1}(\eta_2-\eta_3)=0.
\end{equation}
This implies that ${m_3=0}$. Indeed, if ${m_3>0}$ then~\eqref{esumofthree} can be rewritten as 
\begin{equation}
\label{eimpossible}
\eps_2^{-1}\eta_3^{-1}+\eta_1^{-1}\eps_3^{-1}=-A\eps_1^{-1}(\eta_2-\eta_3).
\end{equation}
Both sides in~\eqref{eimpossible} are non-zero by~\eqref{ezetasdif}, and Lemma~\ref{lsumroots} implies that ${2\ge|A|}$, a contradiction. 
Thus, we have ${m_3=0}$, which, together with~\eqref{ethreelams} and~\eqref{ealleq} implies that 
{%\small
$$
m_1=2\lambda, \quad m_2=\lambda, \quad n_1=\frac32\lambda, \quad n_2=n_3=\frac12\lambda. 
$$}%
Rescaling, we may assume that ${\lambda=2}$, which gives
$$
m_1=4, \quad m_2=2, \quad m_3=0, \quad n_1=3, \quad n_2=n_3=1. 
$$
Using~\eqref{esrut} and~\eqref{esumofthree}, we obtain
{%\small
\begin{align*}
\DD(q)&=
\begin{vmatrix*}[l]
1&1&1\\
\eps_1^{-1}&\eps_2^{-1}q^{2}&\tilf_3q^{4}\\
\eta_1^{-1}&\eta_2^{-1}q^{2}+A\eta_2q^{4}+B\eta_2^2q^5&\eta_3^{-1}q^{2}+A\eta_3q^{4}+B\eta_3^2q^5
\end{vmatrix*}+O(q^6)\\
&=B\eps_1^{-1}(\eta_2^2-\eta_3^2)q^5+O(q^6),
\end{align*}}%
which gives ${\eta_2=\pm\eta_3}$, contradicting~\eqref{esixth}. 

This proves that 
\begin{equation}
\label{einequstr}
2\lambda=m_1-m_3< n_1+n_2.
\end{equation}

\subsection{Proof of ${m_3>0}$}
In addition to this, we have ${m_3>0}$. Indeed, 
equating to~$0$ the coefficient of $q^{2\lambda}$ in (\ref{eterms}), we obtain
\begin{equation}
\label{esumoftwo}
\eps_2^{-1}\eta_3^{-1}+\eta_1^{-1}\tilf_3=0.
\end{equation}
If ${m_3=0}$ then this gives ${f_3=744-\eps_2^{-1}\eta_3^{-1}\eta_1}$, contradicting Lemma~\ref{ljneroot}. 
This proves that 
\begin{equation}
\label{emthpos}
m_3>0,
\end{equation} 
and~\eqref{esumoftwo} rewrites as 
\begin{equation}
\label{esumoftworoots}
\eps_2^{-1}\eta_3^{-1}=-\eps_3^{-1}\eta_1^{-1}.
\end{equation}

\subsection{Proof of ${m_1+m_3=n_1+n_2<3\lambda}$}
\label{ssmmnn<l}
Our next step is showing that ${m_1+m_3=n_1+n_2<3 \lambda}$. 
Using~\eqref{esrut} and~\eqref{esumoftworoots}, we obtain 
{%\small
\begin{align}
\DD(q)&=
\begin{vmatrix*}[l]
1&1&1\\
\eps_1^{-1}&\eps_2^{-1}q^{\lambda}&\eps_3^{-1}q^{2\lambda}+A\eps_3q^{m_1+m_3}\\
\eta_1^{-1}&\eta_2^{-1}q^{\lambda}+A\eta_2q^{n_1+n_2}&\eta_3^{-1}q^{\lambda}+A\eta_3q^{n_1+n_2}
\end{vmatrix*}+o(q^{m_1+m_3})+o(q^{n_1+n_2})\nonumber\\
\label{elongstupid}
&=A\eps_3\eta_1^{-1}q^{m_1+m_3}+A\eps_1^{-1}(\eta_2-\eta_3)q^{n_1+n_2}-\eps_3^{-1}\eta_2^{-1}q^{3\lambda}+o(q^{m_1+m_3})+o(q^{n_1+n_2}).
\end{align}}%

{\sloppy

\paragraph{\underline{We have ${m_1+m_3\ge n_1+n_2}$}}
Indeed, if ${m_1+m_3<n_1+n_2}$ then we must have ${m_1+m_3=3\lambda}$ and ${A\eps_3\eta_1^{-1}=\eps_3^{-1}\eta_2^{-1}}$, which is impossible because~$A$ is not a root of unity.

\paragraph{\underline{We have ${m_1+m_3\le n_1+n_2}$}}
Similarly, if ${m_1+m_3>n_1+n_2}$ then we must have ${n_1+n_2=3\lambda}$ and ${A\eps_1^{-1}(\eta_2-\eta_3)=\eps_3^{-1}\eta_2^{-1}}$, which is impossible because~$A$ does not divide a root of unity. 

}

\paragraph{\underline{We have ${m_1+m_3= n_1+n_2<3\lambda}$}}
Indeed, if ${m_1+m_3=n_1+n_2>3\lambda}$ then the $q^{3\lambda}$ cannot be eliminated. And if ${m_1+m_3=n_1+n_2=3\lambda}$ then  ${A\eps_3\eta_1^{-1}+A\eps_1^{-1}(\eta_2-\eta_3)=\eps_3^{-1}\eta_2^{-1}}$, which is impossible because~$A$ does not divide a root of unity.

\bigskip

This, we proved that 
\begin{equation}
\label{emmnnlll}
m_1+m_3=n_1+n_2<3 \lambda.
\end{equation}
Since ${n_2=n_1-\lambda}$ and ${m_3=m_1-2\lambda}$   (see~\eqref{ethreelams} and~\eqref{einequstr}), this implies that 
{%\small
\begin{equation}
\label{en=m-halfl}
n_1=m_1-\frac12 \lambda. 
\end{equation}}% 
Also, comparing the coefficients in~\eqref{elongstupid}, we obtain
\begin{equation}
\label{eyet}
\eps_3\eta_1^{-1}+\eps_1^{-1}\eta_2-\eps_1^{-1}\eta_3=0.
\end{equation}

\subsection{Conclusion}%{We have ${m_1+2m_3=3\lambda}$}%{Relation between the levels}
\label{ssconcl}

%We want now to prove that ${m_1+2m_3=3\lambda}$. First of all, 
We are almost done. Let us summarize the relations between the levels we already obtained. 
%Setting ${m=m_1}$, we 
We deduce from~\eqref{ethreelams},~\eqref{emthpos},~\eqref{emmnnlll} and~\eqref{en=m-halfl} the following:
$$
%m_1=m, \quad m_2=m-\lambda, \quad m_3=m-2\lambda, \quad n_1=m-\frac12\lambda, \quad n_2=n_3=m-\frac32\lambda,\quad 2\lambda <m <\frac52\lambda. 
m_2=m_1-\lambda, \quad m_3=m_1-2\lambda, \quad n_1=m_1-\frac12\lambda, \quad n_2=n_3=m_1-\frac32\lambda,\quad 2\lambda <m_1 <\frac52\lambda. 
$$
This implies the following inequalities:
$$
2m_1>m_1+m_2=m_1+m_3+\lambda=n_1+n_2+\lambda>3\lambda, \quad 2n_1>3\lambda, \quad %n_1+3n_2=n_1+3n_3>3\lambda, \quad 
n_1+2n_2>m_1+2m_3. 
$$
It follows that 
{%\small
\begin{align*}
\DD(q)&=
\begin{vmatrix*}[l]
1&1&1\\
\eps_1^{-1}&\eps_2^{-1}q^{\lambda}&\eps_3^{-1}q^{2\lambda}+A\eps_3q^{m_1+m_3}+B\eps_3^2q^{m_1+2m_3}\\
\eta_1^{-1}&\eta_2^{-1}q^{\lambda}+A\eta_2q^{n_1+n_2}&\eta_3^{-1}q^{\lambda}+A\eta_3q^{n_1+n_2}
\end{vmatrix*}+o(q^{m_1+2m_3})+o(q^{3\lambda})\\
&=-\eps_3^{-1}\eta_2^{-1}q^{3\lambda}+B\eps_3^2\eta_1^{-1}q^{m_1+2m_3}+o(q^{m_1+2m_3})+o(q^{3\lambda}).
\end{align*}}%
We obtain ${3\lambda=m_1+2m_3}$ and ${\eps_3^{-1}\eta_2^{-1}=B\eps_3^2\eta_1^{-1}}$. But the last equation is impossible because~$B$ is not a root of unity. This proves that~\eqref{emggng=} is impossible in the case~\eqref{eng=g}.

\section{The Case ${m_2>m_3}$,\quad ${n_3>n_2}$}
\label{smggngl}
In this section we assume that 
\begin{equation}
\label{emggngl}
m_1>m_2>m_3,\quad n_1>n_3>n_2,
\end{equation}
(as usual with $m_1 \geq n_1$) and will, eventually, arrive to a contradiction. This is the nastiest case, and we beg for the reader's patience. 

Relation~\eqref{ediffbis} now becomes ${m_1-m_2=n_1-n_3}$. We set
${m_1-m_2=n_1-n_3=\lambda}$.
Using notation~\eqref{etilfg}, we write
{%\small
\begin{align}
\DD(q)&=
\begin{vmatrix*}[l]
1&1&1\\
\eps_1^{-1}&\eps_2^{-1}q^{\lambda}&\tilf_3q^{m_1-m_3}\\
\eta_1^{-1}&\tilg_2q^{n_1-n_2}&\eta_3^{-1}q^{\lambda}
\end{vmatrix*}+o(q^{n_1})\nonumber\\
\label{efirst}
&= \begin{vmatrix}
\eps_1^{-1}&\eps_2^{-1}\\
\eta_1^{-1}&-\eta_3^{-1}
\end{vmatrix}
%(-\eps_1^{-1}\eta_3^{-1}-\eps_2^{-1}\eta_1^{-1})
q^\lambda+\tilf_3\eta_1^{-1}q^{m_1-m_3}+\eps_1^{-1}\tilg_2q^{n_1-n_2}+\eps_2^{-1}\eta_3^{-1}q^{2\lambda}-\tilf_3\tilg_2q^{m_1-m_3+n_1-n_2}+o(q^{n_1}).
\end{align}}%
Since ${0<\lambda<m_1-m_3,n_1-n_2}$, this implies that
{%\small
\begin{equation}
\label{etrd}
%\eps_1^{-1}\eta_3^{-1}+\eps_2^{-1}\eta_1^{-1}
\begin{vmatrix}
\eps_1^{-1}&\eps_2^{-1}\\
\eta_1^{-1}&-\eta_3^{-1}
\end{vmatrix}=0.
\end{equation}}%

\subsection{Proof of ${m_1-m_3=n_1-n_2}$}

Let us start by proving that 
\begin{equation}
\label{etordu}
m_1-m_3=n_1-n_2.
\end{equation}
Indeed, assume that ${m_1-m_3\ne n_1-n_2}$. Then $q^{n_1-n_2}$ in~\eqref{efirst} can be eliminated only if 
\begin{equation}
\label{ejoopa}
n_1-n_2=2\lambda, \quad \eps_1^{-1}\tilg_2=-\eps_2^{-1}\eta_3^{-1}. 
\end{equation}
This implies also that 
${n_2>0}$. 
Indeed, if ${n_2=0}$ then the second equality in~\eqref{ejoopa} gives ${g_2=744-\eps_1\eps_2^{-1}\eta_3^{-1}}$ contradicting Lemma~\ref{ljneroot}.%. Lemma~\ref{ljnecyc} now implies that ${\eps_1\eps_2^{-1}\eta_3^{-1}=\pm1}$ contradicting Lemma~\ref{ljneroot}.%, which gives ${g_2\in\{743,745\}}$, contradicting Lemma~\ref{ljne744}. 

Using~\eqref{etrd} and~\eqref{ejoopa}, we can now write 
{%\small
\begin{align*}
\DD(q)&=
\begin{vmatrix*}[l]
1&1&1\\
\eps_1^{-1}&\eps_2^{-1}q^{\lambda}&\tilf_3q^{m_1-m_3}\\
\eta_1^{-1}&\eta_2^{-1}q^{2\lambda}+A\eta_2q^{n_1+n_2}&\eta_3^{-1}q^{\lambda}
\end{vmatrix*}+o(q^{m_1})+o(q^{n_1+n_2})\\
&= \tilf_3\eta_1^{-1}q^{m_1-m_3}+A\eps_1^{-1}\eta_2q^{n_1+n_2}+o(q^{m_1-m_3})+o(q^{n_1+n_2}).
\end{align*}}%
Here the term with $q^{m_1-m_3}$ cannot be eliminated by $o(q^{n_1+n_2})$ because then $m_1-m_3>n_1+n_2$ and after elimination $q^{n_1+n_2}$ would still be standing. So
\begin{equation}
\label{ejoppa}
m_1-m_3=n_1+n_2, \quad \tilf_3\eta_1^{-1}=-A\eps_1^{-1}\eta_2. 
\end{equation}
However, the second equality in~\eqref{ejoppa} is impossible. Indeed, if ${m_3>0}$ then it becomes ${\eps_3^{-1}\eta_1^{-1}=-A\eps_1^{-1}\eta_2}$, which is clearly impossible because ${A=196884}$ is not a root of unity. And if ${m_3=0}$ then it becomes ${f_3=744-A\eps_1^{-1}\eta_1\eta_2}$, contradicting Lemma~\ref{ljneroot}.%. Lemma~\ref{ljnecyc} now implies that ${\eps_1\eps_2^{-1}\eta_3^{-1}=\pm1}$, which gives ${g_2\in\{744-196884,744+196884\}}$, contradicting Lemma~\ref{ljne744}.

This proves~\eqref{etordu}. We set 
${m_1-m_3=n_1-n_2=\lambda'}$. Since ${m_1\ge n_1}$ by~\eqref{emonegenone}, we may summarize our present knowledge as follows: 
{%\small
\begin{align*}
&m_1>m_2>m_3; \quad n_1>n_3>n_2;\\
&m_1-m_2=n_1-n_3=\lambda>0; \quad m_1-m_3=n_1-n_2=\lambda'>\lambda; \\
& m_1-n_1=m_2-n_3=m_3-n_2\ge 0. 
\end{align*}}%

\subsection{Proof of ${m_3>0}$}

In this subsection we prove that ${m_3>0}$.  We will assume that ${m_3=0}$ and will arrive to a contradiction.

If ${m_3=0}$ then we have
\begin{equation}
\label{eallwrong}
m_1=n_1=\lambda',\quad m_2=n_3, \quad m_3=n_2=0. 
\end{equation}
Using~\eqref{etrd}, we obtain 
{%\small
\begin{align}
\DD(q)&=
\begin{vmatrix*}[l]
1&1&1\\
\eps_1^{-1}&\eps_2^{-1}q^{\lambda}+A\eps_2q^{m_1+m_2}&\tilf_3q^{m_1}\\
\eta_1^{-1}&\tilg_2q^{m_1}&\eta_3^{-1}q^{\lambda}+A\eta_3q^{m_1+m_2}
\end{vmatrix*}+o(q^{m_1+m_2})\nonumber\\
\label{efirstt}
&= \begin{vmatrix}
\eps_1^{-1}&\tilf_3\\
-\eta_1^{-1}&\tilg_2
\end{vmatrix}q^{m_1}+\eps_2^{-1}\eta_3^{-1}q^{2\lambda}+A\begin{vmatrix}
\eps_1^{-1}&\eps_2\\
\eta_1^{-1}&-\eta_3
\end{vmatrix}q^{m_1+m_2}+o(q^{m_1+m_2}).
\end{align}}%
The term with $q^{m_1+m_2}$ can be eliminated if either 
{%\small
\begin{equation}
\label{estrangemat}
\begin{vmatrix}
\eps_1^{-1}&\eps_2\\
\eta_1^{-1}&-\eta_3
\end{vmatrix}=0,
\end{equation}}%
or ${m_1+m_2=2\lambda}$ and 
{%\small
\begin{equation}
\label{etoostrangemat}
A\begin{vmatrix}
\eps_1^{-1}&\eps_2\\
\eta_1^{-1}&-\eta_3
\end{vmatrix}=-\eps_2^{-1}\eta_3^{-1}.
\end{equation}}%
However,~\eqref{etoostrangemat} is impossible because~$A$ does not divide a root of unity. Hence we have~\eqref{estrangemat}. Together with~\eqref{etrd} this implies that 
\begin{equation}
\label{epmpm}
(\eps_1,\eps_2)=\theta
(\eta_1,-\eta_3), \quad \theta=\pm1. 
\end{equation}

The rest of this subsection splits into three cases depending on the relation between~$m_2$ and~$\lambda$.

\paragraph{\underline{The case ${m_2>\lambda}$}}
In this case  ${m_1>2\lambda}$ and $q^{2\lambda}$ in~\eqref{efirstt} cannot be eliminated. 

\paragraph{\underline{The case ${m_2<\lambda}$}} In this case ${m_1<2\lambda}$, and  $q^{m_1}$ in~\eqref{efirstt} can be eliminated only if ${\eps_1^{-1}\tilg_2+\eta_1^{-1}\tilf_3=0}$, which, combined with~\eqref{epmpm}, gives ${\tilg_2=-\theta\tilf_3}$.  Lemma~\ref{lrootofone} implies  that  ${\theta=-1}$ and ${\tilf_3=\tilg_2}$, that is, ${f_3=g_2}$. Also, since ${\theta=-1}$, we obtain ${\eps_2=\eta_3}$, which, together with ${m_2=n_3}$ (see~\eqref{eallwrong})
implies that  ${f_2=g_3}$. This contradicts Lemma~\ref{lcross}.

\paragraph{\underline{The case ${m_2=\lambda}$}}
In this case ${m_1=2\lambda<m_1+m_2}$ and ${\eps_1^{-1}\tilg_2+\eta_1^{-1}\tilf_3+\eps_2^{-1}\eta_3^{-1}=0}$, which contradicts Lemma~\ref{lcute}.

\bigskip

This completes the proof of impossibility of ${m_3=0}$.

\subsection{Proof of ${n_2>0}$}

Thus, we have ${m_3>0}$. Let us now prove that ${n_2>0}$ as well. Indeed, if ${n_2=0}$ then
\begin{equation}
\label{eallwrongagain}
m_1>n_1=\lambda',\quad m_2>n_3, \quad m_3>n_2=0. 
\end{equation}
Using~\eqref{etrd}, we obtain 
{%\small
\begin{align*}
\DD(q)&=
\begin{vmatrix*}[l]
1&1&1\\
\eps_1^{-1}&\eps_2^{-1}q^{\lambda}&\eps_3^{-1}q^{n_1}\\
\eta_1^{-1}&\tilg_2q^{n_1}&\eta_3^{-1}q^{\lambda}
\end{vmatrix*}+o(q^{n_1})= 
(\eps_1^{-1}\tilg_2+\eps_3^{-1}\eta_1^{-1})q^{n_1}+\eps_2^{-1}\eta_3^{-1}q^{2\lambda}+o(q^{n_1}). 
\end{align*}}%
Now to eliminate~$q^{n_1}$ we need to have one of the following:
{%\small
\begin{align}
\label{eonebad}
\eps_1^{-1}\tilg_2+\eps_3^{-1}\eta_1^{-1}&=0,\\
\label{eotherbad}
\eps_1^{-1}\tilg_2+\eps_3^{-1}\eta_1^{-1}+\eps_2^{-1}\eta_3^{-1}&=0. 
\end{align}}%
However, since ${\tilg_2=g_2-744}$, equation~\eqref{eonebad} contradicts Lemma~\ref{ljneroot}. Furthermore, applying Lemma~\ref{ljnetworoots} to equation~\eqref{eotherbad}, we obtain ${g_2\in \{744,744\pm1,744\pm2\}}$, contradicting Lemma~\ref{ljne744etc}. 

This proves that ${n_2>0}$. Let us summarize our present knowledge as follows: 
{%\small
\begin{align*}
&m_1>m_2>m_3>0; \quad n_1>n_3>n_2>0;\\
&m_1-m_2=n_1-n_3=\lambda>0; \quad m_1-m_3=n_1-n_2=\lambda'>\lambda; \\
& m_1-n_1=m_2-n_3=m_3-n_2\ge 0. 
\end{align*}}%

\subsection{Proof of ${m_1=n_1}$}

Next, we show that ${m_1=n_1}$. 
Thus, assume that ${m_1>n_1}$. Then we also have ${m_2>n_3}$ and ${m_3>n_2}$. 
Using~\eqref{etrd}, we write
{%\small
\begin{align}
\DD(q)&=
\begin{vmatrix*}[l]
1&1&1\\
\eps_1^{-1}&\eps_2^{-1}q^{\lambda}&\eps_3^{-1}q^{\lambda'}\\
\eta_1^{-1}&\eta_2^{-1}q^{\lambda'}+A\eta_2q^{n_1+n_2}&\eta_3^{-1}q^{\lambda}
\end{vmatrix*}+o(q^{n_1+n_2})\nonumber\\
%\label{efirst}
&= \begin{vmatrix}
\eps_1^{-1}&\eps_3^{-1}\\
-\eta_1^{-1}&\eta_2^{-1}
\end{vmatrix}q^{\lambda'}+\eps_2^{-1}\eta_3^{-1}q^{2\lambda}-\eps_3^{-1}\eta_2^{-1}q^{2\lambda'}+A\eps_1^{-1}\eta_2q^{n_1+n_2}+o(q^{n_1+n_2}).
\end{align}}%
To eliminate  $q^{n_1+n_2}$ we need one of the following to hold:
{%\small
\begin{align}
\label{elampp}
2\lambda=n_1+n_2, \quad &\eps_2^{-1}\eta_3^{-1}=-A\eps_1^{-1}\eta_2,\\
\label{elamppp}
2\lambda'=n_1+n_2, \quad &\eps_3^{-1}\eta_2^{-1}=A\eps_1^{-1}\eta_2.
\end{align}}%
However, the second equation in both~\eqref{elampp} and~\eqref{elamppp} cannot be true, because~$A$ is not a root of unity.

This proves that ${m_1=n_1}$. Moreover:
\begin{align}
\label{ealloh}
&m_1=n_1>m_2=n_3>m_3=n_2>0; \\
&m_1-m_2=n_1-n_3=\lambda>0; \quad m_1-m_3=n_1-n_2=\lambda'>\lambda. \nonumber
\end{align}

\subsection{Proof of ${\lambda'=2\lambda}$}
Our next quest is proving that ${\lambda'=2\lambda}$.
Using~\eqref{etrd} and~\eqref{ealloh}, we obtain
{%\small
$$
\DD(q)=
\begin{vmatrix*}[l]
1&1&1\\
\eps_1^{-1}&\eps_2^{-1}q^\lambda&\eps_3^{-1}q^{\lambda'}\\
\eta_1^{-1}&\eta_2^{-1}q^{\lambda'}&\eta_3^{-1}q^{\lambda}
\end{vmatrix*}+o(q^{m_1})
=
-\begin{vmatrix}
\eps_1^{-1}&\eps_3^{-1}\\
\eta_1^{-1}&-\eta_2^{-1}
\end{vmatrix}q^{\lambda'}+\eps_2^{-1}\eta_3^{-1}q^{2\lambda}+o(q^{\lambda'}).
$$}%
This already implies that ${\lambda'\le 2\lambda}$; otherwise~$q^{2\lambda}$ cannot be eliminated.

The proof of the opposite inequality ${\lambda'\ge 2\lambda}$ is much more involved. Thus, assume that ${\lambda'<2\lambda}$. Then we must have 
{%\small
$$
\begin{vmatrix}
\eps_1^{-1}&\eps_3^{-1}\\
\eta_1^{-1}&-\eta_2^{-1}
\end{vmatrix}=0.
$$}%
Together with~\eqref{etrd} this implies that 
\begin{equation}
\label{ethecr}
(\eta_1, -\eta_3, -\eta_2)=\theta(\eps_1,\eps_2,\eps_3),
\end{equation}
where~$\theta$ is some root of unity. We obtain 
{%\small
\begin{align*}
\DD(q)&=
\begin{vmatrix*}[l]
\hphantom{\theta^{-1}}1&\hphantom{-\theta^{-1}}1&\hphantom{-\theta^{-1}}1\\
\hphantom{\theta^{-1}}\eps_1^{-1}&\hphantom{-\theta^{-1}}\eps_2^{-1}q^\lambda&\hphantom{-\theta^{-1}}\eps_3^{-1}q^{\lambda'}+A\eps_3q^{m_1+m_3}\\
\theta^{-1}\eps_1^{-1}&-\theta^{-1}\eps_3^{-1}q^{\lambda'}-A\theta\eps_3q^{m_1+m_3}&-\theta^{-1}\eps_2^{-1}q^{\lambda}
\end{vmatrix*}+o(q^{m_1+m_3})\\
&=
-\theta^{-1}\eps_2^{-2}q^{2\lambda}+\theta^{-1}\eps_3^{-2}q^{2\lambda'}+A\eps_3\eps_1^{-1}(\theta^{-1}-\theta)q^{m_1+m_3}+o(q^{m_1+m_3}).
\end{align*}}%
To eliminate $q^{m_1+m_3}$ one of the following should be satisfied:
$$
A\eps_3\eps_1^{-1}(\theta^{-1}-\theta)=\theta^{-1}\eps_2^{-2};\quad A\eps_3\eps_1^{-1}(\theta^{-1}-\theta)=-\theta^{-1}\eps_3^{-2};\quad A\eps_3\eps_1^{-1}(\theta^{-1}-\theta)=0.
$$
Since~$A$ does not divide a root of unity, only the third equation is possible, which implies 
${\theta=\pm1}$. 
If ${\theta=-1}$ then~\eqref{ealloh} and~\eqref{ethecr} imply that ${f_2=g_3}$ and ${f_3=g_2}$, contradicting Lemma~\ref{lcross}. Thus, ${\theta=1}$ and we have
$$
(\eta_1, -\eta_3, -\eta_2)=(\eps_1,\eps_2,\eps_3),
$$
which gives us the following relations:
{%\small
\begin{align*}
q^{m_1}(g_1-744)&=q^{m_1}(f_1-744); \\
q^{m_1}(g_3-744)&=-q^{m_1}(f_2-744)+O(q^{m_1+2m_2}); \\
q^{m_1}(g_2-744)&=-q^{m_1}(f_3-744)+2B\eps_3^2q^{m_1+2m_3}+o(q^{m_1+2m_3}). 
\end{align*}}%
Using this, and the identity 
{%\small
$$
\begin{vmatrix*}[r]
1&1&1\\
a&b&c\\
a&-c+x&-b
\end{vmatrix*}= c^2-b^2+x(a-c),
$$}%
we obtain
{%\small
\begin{align*}
\DD(q)&=
\begin{vmatrix*}[l]
1&\hphantom{-}1&\hphantom{-}1\\
q^{m_1}(f_1-744)&\hphantom{-}q^{m_1}(f_2-744)&\hphantom{-}q^{m_1}(f_3-744)\\
q^{m_1}(f_1-744)&-q^{m_1}(f_3-744)+2B\eps_3^2q^{m_1+2m_3}&-q^{m_1}(f_2-744)
\end{vmatrix*}+o(q^{m_1+2m_3})\\
&=2B\eps_1^{-1}\eps_3^2q^{m_1+2m_3}+ (\eps_3^{-1}q^{m_1-m_3}+A\eps_3q^{m_1+m_3})^2-(\eps_2^{-1}q^{m_1-m_2}+A\eps_2q^{m_1+m_2})^2+o(q^{m_1+2m_3})\\
&=-\eps_2^{-2}q^{2\lambda}+\eps_3^{-2}q^{2\lambda'}+2B\eps_1^{-1}\eps_3^2q^{m_1+2m_3}+o(q^{m_1+2m_3})
\end{align*}}%
(recall that ${\lambda=m_1-m_2}$ and ${\lambda'=m_1-m_3}$). We see that to eliminate $q^{m_1+2m_3}$ we need to have either ${2B\eps_1^{-1}\eps_3^2=\eps_2^{-2}}$ or ${2B\eps_1^{-1}\eps_3^2=-\eps_3^{-2}}$; both are clearly impossible. 

This proves that ${\lambda'=2\lambda}$. Thus, we have 
\begin{equation}
\label{eallho}
m_1=n_1; \quad m_2=n_3=m_1-\lambda; \quad m_3=n_2=m_1-2\lambda>0. 
\end{equation}

\subsection{Proof of ${2\lambda<m_1<3\lambda}$}
\label{ssstwomthree}
%\subsubsection{Conclusion}
Now it is not difficult to show that
\begin{equation}
\label{etwomthree}
2\lambda<m_1<3\lambda. 
\end{equation} 
In fact, ${m_1>2\lambda}$ is already in~\eqref{eallho}.  Next, using~\eqref{etrd}, we obtain
{%\small
\begin{align*}
\DD(q)&=
\begin{vmatrix*}[l]
1&1&1\\
\eps_1^{-1}&\eps_2^{-1}q^\lambda&\eps_3^{-1}q^{2\lambda}+A\eps_3q^{m_1+m_3}\\
\eta_1^{-1}&\eta_2^{-1}q^{2\lambda}+A\eta_2q^{m_1+m_3}&\eta_3^{-1}q^\lambda
\end{vmatrix*}+o(q^{m_1+m_3})\\
&=
(\eps_1^{-1}\eta_2^{-1}+\eps_3^{-1}\eta_1^{-1}+\eps_2^{-1}\eta_3^{-1})q^{2\lambda}
-\eps_3^{-1}\eta_2^{-1}q^{4\lambda}- A
\begin{vmatrix}
\eps_1^{-1}&\eps_3\\
\eta_1^{-1}&-\eta_2
\end{vmatrix}q^{m_1+m_3}+o(q^{m_1+m_3}). 
\end{align*}}%
Since ${m_1>2\lambda}$, this gives 
\begin{equation}
\label{etrd3}
\eps_1^{-1}\eta_2^{-1}+\eps_3^{-1}\eta_1^{-1}+\eps_2^{-1}\eta_3^{-1}=0.
\end{equation}
Further,  if  ${4\lambda <m_1+m_3}$ then $q^{4\lambda}$ cannot be eliminated.  And if  ${4\lambda =m_1+m_3}$ then
{%\small
$$
-\eps_3^{-1}\eta_2^{-1}=  A
\begin{vmatrix}
\eps_1^{-1}&\eps_3\\
\eta_1^{-1}&-\eta_2
\end{vmatrix},
$$}%
which is impossible because~$A$ does not divide a root of unity. 

Thus, we have ${4\lambda>m_1+m_3=2m_1-2\lambda}$, that is, ${m_1<3\lambda}$, proving~\eqref{etwomthree}. 
In addition to this, to eliminate $q^{m_1+m_3}$ we need to have 
{%\small
$$
\begin{vmatrix}
\eps_1^{-1}&\eps_3\\
\eta_1^{-1}&-\eta_2
\end{vmatrix}=0.
$$}%
Together with~\eqref{etrd} this implies that 
\begin{equation}
\label{elasthe}
(\eta_1^{-1},-\eta_3^{-1},-\eta_2)=\theta(\eps_1^{-1},\eps_2^{-1}, \eps_3)
\end{equation}
for some root of unity~$\theta$.

\subsection{Conclusion}
It follows from~\eqref{etwomthree} that ${m_3<\lambda}$, whence 
$$
m_1+2m_3<m_1+m_3+\lambda=m_1+m_2<2m_1. 
$$
Using this, \eqref{etrd},~\eqref{etrd3} and~\eqref{elasthe}, we obtain
{%\small
\begin{align*}
\DD(q)&=
\begin{vmatrix*}[l]
1&1&1\\
\eps_1^{-1}&\eps_2^{-1}q^\lambda&\eps_3^{-1}q^{2\lambda}+A\eps_3q^{m_1+m_3}+B\eps_3^2q^{m_1+2m_3}\\
\eta_1^{-1}&\eta_2^{-1}q^{2\lambda}+A\eta_2q^{m_1+m_3}+B\eta_2^2q^{m_1+2m_3}&\eta_3^{-1}q^\lambda
\end{vmatrix*}+o(q^{m_1+2m_3})\\
&=
-\eps_3^{-1}\eta_2^{-1}q^{4\lambda}- B
\begin{vmatrix}
\eps_1^{-1}&\eps_3^2\\
\eta_1^{-1}&-\eta_2^2
\end{vmatrix}q^{m_1+2m_3}+o(q^{m_1+2m_3}). 
\end{align*}}%
Arguing as in Subsection~\ref{ssstwomthree}, we obtain from this ${4\lambda>m_1+2m_3}$ and
{%\small
\begin{equation*}
\begin{vmatrix}
\eps_1^{-1}&\eps_3^2\\
\eta_1^{-1}&-\eta_2^2
\end{vmatrix}=0,
\end{equation*}}%
which, together with~\eqref{elasthe}, implies that ${\theta=-1}$. It follows that ${\eta_2=\eps_3}$ and ${\eta_3=\eps_2}$; together with~\eqref{ealloh} this implies ${g_2=f_3}$ and ${g_3=f_2}$, contradicting Lemma~\ref{lcross}. 

This completes the proof of impossibility of~\eqref{emggngl}. The Main Lemma is now fully  proved.

\end{document}